\documentclass[11pt,a4paper,rorate,times]{article}
\usepackage{amsfonts}
\usepackage[T1]{fontenc} 
\usepackage[figuresright]{rotating}
\usepackage[bookmarks=true]{hyperref} 
\usepackage{epstopdf}
\usepackage[english]{babel}
\usepackage{latexsym}
\usepackage{amsmath}
\usepackage{amssymb}
\usepackage{algorithm}
\usepackage{algpseudocode}
\usepackage{supertabular}
\usepackage{setspace}
\usepackage{xspace}
\usepackage{url}
\usepackage{color}
\usepackage[table]{xcolor} 
\usepackage{graphicx}
\usepackage{listings}
\usepackage{subcaption}
\usepackage[noblocks]{authblk} 
\usepackage{amsthm}
\usepackage{placeins}
\usepackage{booktabs} 
\usepackage{tabularx} 
\usepackage{ltablex} 
\usepackage{braket}
\usepackage{physics}
\usepackage{threeparttable}

\usepackage[backend = biber, style=authoryear,ibidtracker=context,bibencoding=utf8,hyperref=false,doi=false,url=false, maxcitenames=2, maxbibnames=50, isbn = false, firstinits = true, uniquename=false, uniquelist=false,dashed=false,natbib=true]{biblatex}
\renewbibmacro{in:}{}
\AtEveryBibitem{%
	\clearfield{note}%
}
\usepackage{ifpdf}
\normalfont
\makeatletter

\usepackage{algorithm,algcompatible}
\algnewcommand\INPUT{\item[\textbf{Input:}]}%
\algnewcommand\OUTPUT{\item[\textbf{Output:}]}%

\usepackage{algpseudocode} 
\usepackage{blkarray}  
\usepackage{arydshln}  
\usepackage{mathtools} 
\newcommand{\breakingcomma}{
	\begingroup\lccode`~=`,
	\lowercase{\endgroup\expandafter\def\expandafter~\expandafter{~\penalty0 }}}
\usepackage{breqn}     

\usepackage[toc,page]{appendix}
\usepackage{multirow}
\usepackage{makecell}
\usepackage{colortbl}
\usepackage{adjustbox}
\usepackage{scrextend}

\definecolor{gray75}{gray}{0.75}
\definecolor{gray70}{gray}{0.70}
\definecolor{gray65}{gray}{0.65}
\definecolor{gray60}{gray}{0.60}
\definecolor{gray55}{gray}{0.55}
\definecolor{gray50}{gray}{0.50}
\definecolor{gray45}{gray}{0.45}
\definecolor{gray40}{gray}{0.40}
\definecolor{gray35}{gray}{0.35}
\definecolor{gray30}{gray}{0.30}

\definecolor{purple}{RGB}{160, 32, 240} 
\definecolor{magenta}{RGB}{255, 0, 255} 
\definecolor{orange}{RGB}{255, 165, 0}

\renewcommand{\thefootnote}{\arabic{footnote}}

\newcommand{\astfootnote}[1]{
	\let\oldthefootnote=\thefootnote
	\setcounter{footnote}{0}
	\renewcommand{\thefootnote}{\fnsymbol{footnote}}
	\footnote{#1}
	\let\thefootnote=\oldthefootnote
	\setcounter{footnote}{4}
}

\begin{document}
	\thispagestyle{empty}
			\title{Benchmarking of Quantum and Classical Computing in Large-Scale Dynamic Portfolio Optimization Under Market Frictions}
\author[1]{Ying Chen}
\author[2]{Thorsten Koch}
\author[3]{Hanqiu Peng}
\author[3]{Hongrui Zhang}
\affil[1]{Centre for Quantitative Finance, Department of Mathematics \& Asian Institute of Digital Finance, \& Risk Management Institute, National University of Singapore}
\affil[2]{Zuse Institute Berlin \& Technische Universit\"at Berlin}
\affil[3]{Department of Mathematics, National University of Singapore}

\maketitle

\date{}

\begin{abstract}
Quantum computing is poised to transform the financial industry, yet its advantages over traditional methods have not been evidenced. As this technology rapidly evolves, benchmarking is essential to fairly evaluate and compare different computational strategies. 
This study presents a challenging yet solvable problem of large-scale dynamic portfolio optimization under realistic market conditions with frictions. We frame this issue as a Quadratic Unconstrained Binary Optimization (QUBO) problem, compatible with digital computing and ready for quantum computing, to establish a reliable benchmark. By applying the latest solvers to real data, we release benchmarks that help verify true advancements in dynamic trading strategies, either quantum or digital computing, ensuring that reported improvements in portfolio optimization are based on robust, transparent, and comparable metrics.

\end{abstract}


\newpage
\section{Introduction}
Quantum computing is depicted as a remarkable next-generation technology that could lead to significant breakthroughs in the financial industry, promising to transform traditional practices. However, its definitive advantages over classical computing methods have not yet been empirically proven. As this technology continues to develop, the role of benchmarking becomes increasingly important. Effective benchmarking is essential for verifying the real benefits of new computational technologies, ensuring that their advantages are not merely due to comparisons with outdated methods. This work presents a challenging yet solvable problem in the realm of dynamic trading strategies, where large-scale portfolio optimizations contend with realistic market conditions and frictions known to impact performance. We formulate the problem within the framework of Quadratic Unconstrained Binary Optimization (QUBO), solvable with both quantum and digital computing. This allows that any reported advantage is measurable within both computing time and precisely defined optimality bounds, providing a standardized approach for measuring and comparing the efficacy of computational strategies on a transparent, equivalent, and robust basis.

Dynamic trading strategies are designed to continuously adjust asset allocation in response to market changes, aiming to capitalize on opportunities in the ever-changing financial markets. The importance of dynamic trading strategies lies in their ability to adapt swiftly to market dynamics, thus potentially enhancing portfolio performance and mitigating risks more effectively than their static counterparts \citep{fama1970efficient, sharpe1998sharpe, brown2022information}. These strategies are inherently active and efficient, making them vital in investment and trading decision, particularly when dealing with a diverse array of stocks \citep{gatzert2021portfolio, herman2022survey}. However, from a computing perspective, the practical implementation of these strategies presents significant challenges. The complexities involved include integrating transaction costs, managing constraints like short selling, and executing large-scale optimizations that involve numerous parameters and a large number of stocks \citep{pliska1997introduction,campbell1998econometrics,mansini2015linear}. Additionally, it often needs to accommodate integer variables for specifics such as the number of stocks to be selected, share counts or trading blocks. These technical requirements add layers of complexity to the computational algorithms needed, making them computationally intensive and challenging. This not only increases the demand on computational resources but also requires advanced optimization techniques to solve these problems effectively within reasonable time frames, see e.g. \textcite{bertsimas1997introduction}, \textcite{lin2008genetic}, \textcite{liu2015multi}, and \textcite{dai2021learning}.

Continuous models in portfolio optimization have gained popularity in the literature since the introduction of Markowitz's Modern Portfolio Theory in 1952. Originating from single-period stochastic control \citep{samuelson1975lifetime,merton1975optimum,dybvig1985differential} and evolving to address multi-period settings \citep{li2000optimal, brandt2006dynamic}, these models assume unlimited transaction frequency and provide critical insights through mathematical derivations. In certain idealized settings, they can offer optimal trading strategies with analytical solutions that facilitate straightforward plug-in implementation. However, the practical applicability of continuous models faces several limitations. Firstly, actual trading operates under discrete conditions, with limited trading hours, sizes, and frequencies, which reduce the practical relevance of continuous models in real-world scenarios. Additionally, when factors such as transaction costs, taxation, or heavy-tailed distributions are incorporated, continuous models often struggle to maintain their analytical tractability. Adapting these models to discrete settings can lead to loss of precision, heightened computational demands, and potential instability. The numerical complexity typically constrains their use in portfolio analysis to very few assets.

Discrete models in portfolio optimization offer a more realistic representation of financial markets by closely mirroring actual trading practices, including the granularity required for specific constraints, such as integer constraints on the number of shares \citep{bertsimas1998optimal,brown2011dynamic}. These models align with real-world operations, providing the flexibility to handle market frictions such as transaction costs and short selling. Nevertheless, this increased realism introduces significant computational challenges. When managing portfolios with a large number of stocks, compounded by market frictions, the optimization problem becomes NP-hard, characterized by slow convergence rates and high time complexity \citep{broadie1998optimal,zabarankin2014capital}. In practice, the number of stocks is often limited to around two digits. As new financial products continuously emerge and trading frequencies increase, the disparity between theoretical models and practical implementation underscores the urgent need for technologies that can efficiently solve dynamic trading strategies for large-scale portfolio optimization problems.

In 2012, Scott introduced the term ``quantum supremacy,'' which spurred subsequent research showcasing the theoretical potential of utilizing quantum computing \citep{nielsen2010quantum, harrow2017quantum, mcgeoch2022adiabatic}. This has attracted attention across various fields, including credit rating and derivatives pricing \citep{rebentrost2018quantum, woerner2019quantum, orus2019quantum, egger2020credit, stamatopoulos2020option}. Theoretical results also indicate a potential polynomial speedup in portfolio optimization under certain conditions, such as bounded condition number for Newton matrices and tomography precision limited by the duality gap \citep{kerenidis2019quantum}. Table \ref{tab:algo_sum} summarizes works involving experimental computing on dynamic trading strategies for portfolio optimization using both digital and quantum computing. Current quantum studies focus on discrete optimization with integer decision variables but frequently overlook transaction costs, prohibit short selling, and limit the number of assets to around 10. This suggests that, although quantum portfolio optimization has potential, it is still in its early stages \citep{elsokkary2017financial, braine2021quantum, gilliam2021grover, rosenberg2016solving, mugel2022dynamic}. 

In fact, solving large-scale dynamic portfolio optimization problems is also challenging for classical computing under realistic settings. This scenario highlights the appropriateness of benchmarking computational methods and algorithms for large-scale optimization problems: It presents a significant challenge; the potential for solution exists, and the problems, while complex, are solvable and importance. A critical aspect is determining the proper settings for the benchmark problem. We consider a multi-period discrete portfolio optimization problem that incorporates factors such as short sales, transaction costs, integer variables, and large-scale optimization challenges. We formulate this problem in the form of Quadratic Unconstrained Binary Optimization (QUBO), where the control variable is a vector of binary variables, making any 0/1 vector a feasible solution. QUBO is a well-established field in optimization where digital algorithms have been developed to reliable standards. This allows not only for obtaining a solution but also for the precise measurement of its optimality. Currently, we manage portfolios of up to 499 stocks, with each stock selectable up to three times over 15 time steps. This configuration is challenging yet manageable and serves to delineate the current state-of-the-art in digital computing.

Furthermore, there is a direct link between QUBO and quantum computing. QUBO problems can be transformed into Ising Hamiltonians, enabling the direct implementation of quantum algorithms using QUBO inputs suitable for noisy intermediate-scale quantum (NISQ) hardware through variational quantum simulation (VQS) methods such as the Variational Quantum Eigensolver (VQE) and the Quantum Approximate Optimization Algorithm (QAOA). For instance, QAOA is explicitly designed to solve QUBOs, and theoretically, running sufficient steps should yield a high probability of obtaining a good or optimal solution. Although current quantum computing settings do not yet effectively solve these problems, the rapid advancement of quantum computing hardware and software may soon yield a quantum advantage. This benchmarking aims to accurately ascertain the genuine state-of-the-art in both digital and quantum computing, providing a measurable and comparable basis. Understanding these advancements is crucial as it prepares us for future developments and allows us to leverage potential benefits as quantum hardware continues to evolve.

Our contribution focuses on developing a benchmarking framework to evaluate dynamic trading strategies within a multi-period discrete portfolio optimization problem under realistic market conditions. We aim to provide fair and consistent measurement across computational technologies, including both digital and quantum computing. By modeling the problem as a QUBO, which is challenging yet solvable with current computing technologies, we enable rigorous comparisons across various computational methods. Additionally, we assess the numerical performance of the latest optimization solvers using real-world data to establish benchmarks for large-scale dynamic portfolio optimization. This work seeks to bridge the gap between the theoretical potential of advanced computing technologies and their practical applicability in financial markets, emphasizing the importance of benchmarking as computational technologies continue to evolve.

\begin{table}[ht]
\centering
{\footnotesize\begin{tabular}{l|c|c|c|c|c|c|c}
\hline\hline
 & \textbf{\#}& \textbf{Trac.} & \textbf{Max} & \textbf{Time} & \textbf{Short} & \textbf{Int.} &  \\ 
\textbf{Authors}  & \textbf{Stocks} & \textbf{Cost} & \textbf{Shares} & \textbf{Steps} & \textbf{Sell} & \textbf{Var} & \textbf{Algo}\\
\hline
Lobo et al., 2007 (C) & 100 & {\color{green}Yes} & - & 1 & {\color{green}Yes} & {\color{red}No} & Convex opti. \\
Wang and Zhou, 2020 (C) & 1 & {\color{red}No} & - & - & {\color{red}No} & {\color{red}No} & RL \\
Muthuraman et al., 2008 (C) & 2 & {\color{green}Yes} & - & - & {\color{red}No} & {\color{red}No} & DP\\
Dai et al., 2010 (C) & 1 & {\color{green}Yes} & - & - & {\color{red}No} & {\color{red}No} & DP\\ \hline
Bertsimas and Lo, 1998 (D) & 1 & {\color{green}Yes} & - & 20 & {\color{red}No} & {\color{green}Yes} & DP \\
Guastaroba et al., 2020 (D) & 100 & {\color{red}No} & - & 104 & {\color{red}No} & {\color{red}No} & LP \\
Zabarankin et al., 2014 (D) & - & {\color{red}No} & - & 90 & {\color{red}No} & {\color{red}No} & LP \\
Brown et al., 2011 (D) & 10 & {\color{green}Yes} & - & 48 & {\color{red}No} & {\color{red}No} & DP \\ \hline
Mugel et al. 2022 (Q) &  8 & {\color{green}Yes} & 7 & 53 & {\color{red}No} & {\color{green}Yes} & Hybrid, VQE, TN \\ 
Brandhofer et al. 2022 (Q) & 10 & {\color{red}No} & 1 & 5 & {\color{red}No} & {\color{green}Yes} & QAOA \\ 
Rosenberg et al. 2016 (Q) & 6 & {\color{red}No} & 3 & 6 & {\color{red}No} & {\color{green}Yes} & Annealing \\ 
Slate et al. 2021 (Q) & 8 & {\color{red}No} & 1 & 1 & {\color{red}No} & {\color{green}Yes} & QWOA \\ \hline
Benchmark & 499 & {\color{green}Yes} & 3 & 15 & {\color{green}Yes} & {\color{green}Yes} & BQP, QUBO\\\hline\hline
\end{tabular}}
\caption{Summary of Experimental Computing on Portfolio Optimization. Digital computing (continuous (C) and discrete (D) models) and quantum (Q) computing. Note: RL: Reinforcement Learning. DP: Dynamic Programming. LP: Linear Programming.}
\label{tab:algo_sum}
\end{table}

The remainder of this paper is organized as follows. Section \ref{sec:model} presents the model setting and QUBO framework. Section \ref{sec:quantum} presents the quantum algorithm. Section \ref{sec:emp} showcases the numerical results and the empirical findings. Section \ref{sec:conc} concludes.
 
\section{Model Settings}\label{sec:model}
In this section, we present the mathematical framework for a portfolio optimization problem over a finite investment horizon divided into multiple time periods. The model includes \(n\) risky assets alongside one risk-free asset (bank deposit). At each period, the investor has the opportunity to reallocate assets within the portfolio, which incurs transaction costs. Short selling is permitted. The objective of the investor is to maximize total wealth and minimize risk across the investment horizon, while satisfying several constraints: the availability of capital throughout the period, a ceiling on the number of distinct assets that can be held, and limits on the quantity of each asset. The formulation aims to develop a discrete-time dynamic portfolio allocation strategy that conforms to the Markowitz mean-variance objective function.

\subsection{Classic Mean-Variance Portfolio Optimization Problem}
In the classic Markowitz's portfolio optimization problem, the simplest form requires choosing exactly \(B\) assets out of \(n\), which can be formulated as finding the binary decision variables \(x_i \in \{0,1\}\) for each asset \(i\) that maximize the quadratic utility over one single period:
\begin{equation*}
\max_{\mathbf{x} \in \{0,1\}^n} \boldsymbol{\mu}^{\mathsf{T}} \mathbf{x} - q \mathbf{x}^{\mathsf{T}} \boldsymbol{\Sigma} \mathbf{x} \quad \text{subject to} \quad \mathbf{1}^{\mathsf{T}} \mathbf{x} = B,
\end{equation*}
where \( \boldsymbol{\mu} \) represents the vector of expected returns, \( \boldsymbol{\Sigma} \) is the covariance matrix of asset prices, and \( q > 0\) is the so-called risk tolerance factor. This factor balances the risk against expected returns in the portfolio; a higher value of \(q\) emphasizes risk aversion, leading to a more conservative investment strategy, while a value of \(q=0\) completely ignores risks. 

Using a sufficiently large penalty \(P\) and reversing the objective, the problem can be reformulated as an unconstrained one:
\[
\min_{\mathbf{x} \in \{0,1\}^n} -\boldsymbol{\mu}^{\mathsf{T}} \mathbf{x} + q \mathbf{x}^{\mathsf{T}} \boldsymbol{\Sigma} \mathbf{x} + P(B - \mathbf{1}^{\mathsf{T}} \mathbf{x})^2,
\]
or equivalently in element-wise form:
\[
\min_{x_i \in \{0,1\}} -\sum_i \mu_i x_i + q \sum_i \sum_j x_i \sigma_{ij} x_j + P \left( B - \sum_i x_i \right)^2.
\]
Note that this single-period problem takes the form of Quadratic Unconstrained Binary Optimization (QUBO), where the objective is to minimize a quadratic polynomial 
\[\min_{\mathbf{x} \in \{0,1\}^n} \mathbf{x}^{\mathsf{T}} \mathbf{Q} \mathbf{x} = \sum_{i=1}^n \sum_{j=1}^n Q_{ij} x_i x_j,\]
where \(\mathbf{Q}\) is a matrix that is dense and grows quadratically with the number of assets \(n\). 

\subsection{Multiple Time Periods with Transaction Costs and Short Selling}
Consider an investment horizon comprising multiple time periods \(t \in \{1, \dots, T\}\), each separated by a unit length of time. Introduce an initial period \(0\) with the asset selection \(x_{i0} = 0\). The investor is to allocate the available (normalized) units of capital \(C\) at each period, by either reallocating up to \(B\) assets or holding cash, subject to the constraint \(C \leq B\). Let \(p_{i,t}\) represent the price of one unit of asset \(i\) at time \(t\), and let \(\sigma_{ij,t}\) denote the covariance between stocks \(i\) and \(j\) at time \(t\). The portfolio risk for this model is defined by:
\begin{equation}
\sum_i \sum_j p_{it} x_{it} \sigma_{ijt} x_{jt} p_{jt},
\end{equation}
The task is to determine the trajectory \(x_{it} \in \{0,1\}\) that minimizes the overall portfolio risk while maximizing profit, expressed as:
\begin{equation}
\sum_{t=1}^T \left( q \sum_i \sum_j p_{it} x_{it} \sigma_{ijt} x_{jt} p_{jt} - \sum_i (p_{i,t+1} - p_{it}) x_{it} \right),
\end{equation}
We allow that each asset \(i\) can be selected up to \(k\) times and define a \(2kn\)-dimensional binary vector \(\mathbf{x}\) for the \(n\) assets at each time \(t\), where each block of \(k\) binary variables represents a single asset type. The inclusion of a factor of 2 accounts for the allowance of short selling, with an equal number of vector elements allocated to represent both long and short positions. To mirror realistic trading scenarios, constraints on transaction costs and short selling are imposed at each time step. Both impose limits on trading activities, and affect the path-dependent cash.

\noindent\textbf{Transaction Costs}: We denote \(\delta\) as the transaction cost rate applied to both buying and selling. Since \(x_{it}\) is binary, its absolute value change is given by \(\lvert x_{it-1} - x_{it} \rvert = x_{it-1} + x_{it} - 2 x_{it-1} x_{it}\). The total transaction cost is represented by:
\begin{equation}
\delta \sum_i \sum_{t=1}^T p_{it} (x_{it-1} + x_{it} - 2 x_{it-1} x_{it}).
\end{equation}

\noindent We also account for the transaction cost associated with the liquidation of all assets at the terminal time \(T\):
\begin{equation}
\delta \sum_i p_{iT} x_{iT}.
\end{equation}

\noindent\textbf{Short Selling Cost and Restriction}: We introduce a short-selling indicator, \(\tau \in\{-1, +1\}\), where \(-1\) indicates a short position and \(+1\) a long position. Our model includes a borrowing cost rate for short sales, denoted by \(\rho_s\), to enforce prudence in short-selling activities:
\begin{equation}
\rho_s \sum_{i \in S} p_{it} x_{it},
\end{equation}
where S denotes the set of assets selected for short selling. Each short position allows the acquisition of additional capital for investment. We introduce slack variables \(s_{bt} \in \{0,1\}, b \in \{0, \ldots, \lfloor \log_2 B \rfloor\}\) to cap the total assets:
\begin{equation}\label{eq:ss}
\sum_i x_{it} + \sum_b 2^b s_{bt} = B \quad \text{for all } t \in \{1, \ldots, T\}.
\end{equation}

\noindent\textbf{Path-dependent Cash}: In each time period, the model restricts the total available cash to not exceed \(C\) units. We introduce another slack variables \(y_{ct} \in \{0, 1\}\), \(c \in \{0, \ldots, \lfloor \log_2 C \rfloor\}\). The constraint for cash, given short-selling is permitted, is formulated as:
\begin{equation}\label{eq:cash}
\sum_i \tau_i x_{it} + \sum_c 2^c y_{ct} = C \quad \text{for all } t \in \{1, \ldots, T\}.
\end{equation}
Note that the slack variables \(s_{ct}\) indicate the amount of the budget that remains unspent and incur interest:
\[
- \rho_c u \sum_{t} \sum_{c} 2^c y_{ct},
\]
where \(\rho_c\) denotes the risk-free interest rate, and \(u\) is the value of one normalized cash unit.
\noindent The constraints (\ref{eq:ss}) and (\ref{eq:cash}) can be enforced in the objective function through a quadratic penalty term:
\begin{equation}
P \sum_{t=1}^T \left\{\big( B - \sum_i x_{it} - \sum_b 2^b s_{bt} \big)^2 + 
\big( C - \sum_i \tau_i x_{it} - \sum_c 2^c y_{ct} \big)^2 \right\}.
\end{equation}

\subsection{QUBO}
For multi-period portfolio optimization under market frictions, it is not only challenging to obtain a solution when handling large dimensional assets, but also difficult to prove exact global optimal solutions. Quadratic Unconstrained Binary Optimization (QUBO) provides a framework suitable for a broad range of combinatorial optimization problems. The adoption of QUBO in benchmarking this optimization problem is justified for several reasons: Above all, QUBO is a vibrant area in optimization where digital algorithms have been developed to reliable standards, allowing us to not only obtain solutions but also measure their optimality with a precise bound. More importantly, QUBO is quantum-ready. Any QUBO problem can be transformed into Ising Hamiltonians, allowing direct implementation of quantum algorithms using QUBO inputs on quantum computing platforms like IBM Qiskit, D-Wave annealing, or hybrid systems.

The portfolio optimization problem can be formulated in QUBO form as follows:
\begin{eqnarray}
&\min_{\substack{x \in \{0,1\}^{n \times t} \\ y \in \{0,1\}^{c \times t} \\ s \in \{0,1\}^{b \times t}}} \sum_{t=1}^{T} \left(q \underbrace{\sum_i \sum_j p_{it} x_{it} \sigma_{ijt} x_{jt} p_{jt}}_{\text{risk}} - \sum_i \Big(\underbrace{(p_{it+1} - p_{it}) x_{it}}_{\text{profit}} - \underbrace{\delta p_{it} (x_{it-1} + x_{it} - 2 x_{it-1} x_{it})}_{\text{transaction cost}}\Big) \right. \\ \nonumber
& \left. \hspace{2.5cm} - \underbrace{\rho_c u \sum_c 2^c y_{ct}}_{\text{cash interest}} + \underbrace{\rho_s \sum_{i \in S} p_{it} x_{it}}_{\text{short selling cost}}  + \underbrace{\delta \sum_i p_{iT} x_{iT}}_{\text{liquidation cost}} \right) \\\nonumber
& \hspace{2.5cm} + P \sum_{t=1}^{T} \left( \underbrace{\left( C - \sum_i \tau_i x_{it} - \sum_c 2^c y_{ct} \right)^2}_{\text{capital limit}} + \underbrace{\left( B - \sum_i x_{it} - \sum_b 2^b s_{bt} \right)^2}_{\text{number of assets limit}} \right)
\end{eqnarray}

The problem can be formulated as an equivalent Binary Quadratic Program (BQP) where the objective is a quadratic function subject to linear constraints, and the decision variables are constrained to be binary. 
\begin{eqnarray}
&\min_{\substack{x \in \{0,1\}^{n \times t} \\ y \in \{0,1\}^{c \times t} \\ s \in \{0,1\}^{b \times t}}} \sum_{t=1}^{T} \left( q \underbrace{ \sum_i \sum_j p_{it} x_{it} \sigma_{ijt} x_{jt} p_{jt}}_{\text{risk}} - \sum_i \Big(\underbrace{(p_{it+1} - p_{it}) x_{it}}_{\text{profit}} - \underbrace{\delta p_{it} (x_{it-1} + x_{it} - 2 x_{it-1} x_{it})}_{\text{transaction cost}}\Big)\right. \\\nonumber
&\left. - \underbrace{\rho_c u \sum_c 2^c y_{ct}}_{\text{cash interest}} + \underbrace{\rho_s \sum_{i \in S} p_{it} x_{it}}_{\text{short selling cost}} + \underbrace{\delta \sum_i p_{iT} x_{iT}}_{\text{liquidation cost}} \right)
\end{eqnarray}
subject to
\[
\sum_i \tau_i x_{it} + \sum_c 2^c y_{ct} = C \quad \text{for all } t \in \{1, \ldots, T\} \quad \quad {\text{capital limit}}
\]
\[
\sum_i x_{it} + \sum_b 2^b s_{bt} = B \quad \text{for all } t \in \{1, \ldots, T\} \quad \quad {\text{number of assets limit}}
\]
The BQP can only be solved by classical digital computing, but with precise bound and reliable performance. It can serve as a baseline to evaluate the optimal solutios of QUBO via digital and quantum computing.

\section{Digital and Quantum-Ready Algorithms for Optimization}\label{sec:quantum}

As a benchmark, we focus on the effectiveness of optimizing BQP and QUBO problems without prioritizing either quantum or digital computing. If solved to optimality, the objective function value should be consistent, irrespective of the modeling approach or the solver used\footnote{A solution may be optimal, guaranteed (bounded), or simply feasible. In all three cases, the solution may be the same, but our understanding and confirmation of its quality differ. Proving optimality often requires substantial computational effort, but solution times can vary significantly. If the goal is to find a good solution, this can be achieved more efficiently.}.

\subsection{Digital Algorithms}

Binary Quadratic Programming (BQP) problems are challenging to solve due to their NP-hard nature. To tackle these problems, a variety of algorithms are typically employed, including Branch-and-Bound, Cutting Planes, and heuristic methods. The Branch-and-Bound algorithm, in particular, is widely utilized for solving BQP problems by systematically exploring and pruning the solution space. Algorithm \ref{alg:BB_BQP} describes the procedure of the Branch-and-Bound algorithm.
\begin{algorithm}[H]
\caption{Branch-and-Bound for BQP Problems}
\label{alg:BB_BQP}
\begin{algorithmic}[1]
\STATE Initialize the best known solution $x^*$ with a large objective value $z^* \gets \infty$
\STATE Initialize the list of subproblems with the original problem
\WHILE{there are unexplored subproblems}
    \STATE Select and remove a subproblem $P$ from the list
    \STATE Solve the relaxation of $P$ to obtain the lower bound $z_P$
    \IF{$z_P \geq z^*$}
        \STATE Discard subproblem $P$
    \ELSE
        \IF{the solution is feasible and $z_P < z^*$}
            \STATE Update the best known solution $x^* \gets x_P$
            \STATE Update the best known objective value (upper bound) $z^* \gets z_P$
        \ELSE
            \STATE Branch on a decision variable to create new subproblems $P_1, P_2$
            \STATE Add $P_1$ and $P_2$ to the list of subproblems
        \ENDIF
    \ENDIF
\ENDWHILE
\STATE \textbf{return} $x^*$ as the optimal solution
\end{algorithmic}
\end{algorithm}
BQP problems can potentially be solved to optimality using advanced solvers like Gurobi, CPLEX, COPT, and SCIP. In benchmark experiments, Gurobi 11.0 was used on an AMD Ryzen 9 Workstation, providing a bound to indicate the optimality gap. Specifically, the optimality gap measures the difference between the upper and lower bounds on the objective function value. While the upper bound is obtained from a feasible solution where the decision variables are binary integers, the lower bound is derived by relaxing the integrality constraints of the decision variables, allowing them to take continuous values within the range \([0,1]\). This relaxation makes the problem easier to solve but may not yield a feasible solution. The gap serves as an indicator of how close the best-known feasible solution is to the optimal solution of the relaxed problem. If the gap is zero, the feasible binary solution is proven to be optimal; if the gap is positive, there could still be a better binary solution, or it might already be the best possible feasible solution but cannot be proven. 

Adaptive Bulk Search (ABS2), proposed by \textcite{nakano2023diverse}, is a framework that combines multiple search algorithms, genetic operations, and solution pools to explore the solution space of QUBO problems. It dynamically adapts strategies based on solution performance, leveraging GPU power to enhance search efficiency and diversity. Algorithm \ref{alg:ABS2} outlines the steps of ABS2, which has been integrated into the qoqo solver designed by \textcite{rehfeldt2023faster}, and was utilized on an NVIDIA A100-SXM4-80GB GPU in our study. By comparing the objective value to the lower bound of BQP, we can assess whether a QUBO solver achieves or exceeds the performance level of the optimal solution computed by a BQP solver.
\begin{algorithm}[H]
\caption{ABS2 Algorithm for QUBO Problems}
\label{alg:ABS2}
\begin{algorithmic}[1]
\STATE Initialize multiple solution pools, each with a random feasible solution $X$
\FOR{each solution pool on a GPU}
    \STATE Compute initial energy $E(X)$ and set $BEST \gets X$, $E(BEST) \gets E(X)$
    \STATE Initialize the energy change for flipping each bit in $X$
\ENDFOR
\REPEAT
    \FOR{each solution pool on a GPU}
        \FOR{each CUDA block (in parallel)}
            \STATE Select a search algorithm and genetic operation based on past performance
            \STATE Generate and refine candidate solutions using genetic operations and bit-flipping
            \STATE Update the energy $E(X)$ by applying the calculated change for the flipped bit
            \IF{$E(X) < E(BEST)$} \STATE Update $BEST \gets X$, $E(BEST) \gets E(X)$
            \ENDIF
            \STATE Return the best solutions to the host and update the solution pool
        \ENDFOR
    \ENDFOR
\UNTIL{stopping criterion is met (e.g., max iterations, convergence, time limit)}
\STATE \textbf{return} the best solution $x^*$ across all pools
\end{algorithmic}
\end{algorithm}

\subsection{Quantum Algorithms}

Quantum algorithms, such as the Quantum Approximate Optimization Algorithm (QAOA), are tailored to solve QUBOs by converting the problem into an Ising model \citep{lucas2014ising}. In the Ising model, variables take the values of "spin up" (↑) and "spin down" (↓), which correspond to +1 and -1, respectively. The objective function in the Ising model is formulated as:
\[
    F(s)=\sum_{i=1}^n h_i s_i + \sum_{i<j} J_{ij} s_i s_j,
\]
where \( h_i \) are the linear coefficients representing qubit biases and \( J_{ij} \) are the quadratic coefficients representing coupling strengths. The conversion from QUBO to the Ising model uses the simple arithmetic mapping \( x_i = \frac{s_i + 1}{2} \), allowing a direct relationship between the two formulations.  

Variational Quantum Eigensolver (VQE) is a hybrid algorithm that combines quantum and classical resources to find the ground state energy of quantum systems. The VQE process (Algorithm \ref{alg:VQE}) involves initializing a parametrized quantum circuit, executing the circuit on a quantum processor to estimate the energy, using a classical optimizer to adjust the circuit parameters, and repeating the process until the energy converges to a minimum \citep{peruzzo2014variational}. VQE is particularly suited for near-term quantum devices and can be designed for QUBO problems using specific rotation and entanglement gates to balance expressiveness and noise resilience.
\begin{algorithm}[H]
\caption{VQE Algorithm for QUBO Problems}
\label{alg:VQE}
\begin{algorithmic}[1]
\STATE $H \gets \text{QUBO to Hamiltonian}(Q)$
\STATE $n \gets \text{Number of qubits}(Q)$
\STATE $\text{ansatz} \gets \text{Create Variational Circuit}(n)$
\STATE $\theta_0 \gets \text{Initialize Parameters}()$
\STATE $\text{optimizer} \gets \text{Choose Optimizer}()$
    
\STATE \textbf{Objective Function}($\theta$):
    \STATE \quad $|\psi(\theta)\rangle \gets \text{Prepare State}(\text{ansatz}, \theta)$
    \STATE \quad $E(\theta) \gets \langle\psi(\theta)|H|\psi(\theta)\rangle$
    \STATE \quad \textbf{return} $E(\theta)$
    
\STATE $\theta_{\text{opt}}, E_{\text{min}} \gets \text{optimizer.Minimize}(\text{Objective Function}, \theta_0)$
\STATE $x_{\text{QUBO}} \gets \text{Convert to QUBO Solution}(\theta_{\text{opt}})$
\STATE \textbf{return} $E_{\text{min}}, x_{\text{QUBO}}$

\end{algorithmic}
\end{algorithm}
QAOA, one of the well-known VQE algorithm, is specifically designed to solve QUBOs. QAOA involves creating a superposition of all possible solutions using a quantum circuit, applying problem-specific unitary operators \( U(Q, q) = e^{-iq Q} \) and mixer Hamiltonians \( U(H, \beta) = e^{-i\beta H} \), measuring the quantum state to obtain an approximate solution, and adjusting the parameters \( q \) and \( \beta \) with a classical optimizer to iteratively improve the solution. The effectiveness of QAOA is influenced by the choice of classical optimizer, quantum circuit depth, initial optimization points, and mixer Hamiltonian selection. Theoretically, running sufficient steps should yield a high probability of obtaining a good or optimal solution. However, current quantum computers cannot execute many steps. Recent papers, such as \cite{montanez2024towards}, have claimed that it is possible to reduce the number of steps required. Ultimately, the more efficiently this reduction is achieved, the fewer steps are needed to reach a high probability of finding a good solution.

In this paper, VQE algorithm is implemented using IBM Qiskit, which offers access to gate-based superconducting quantum processors. IBM's quantum processors contain over 100 qubits. The qubits are arranged in specific topologies that allow for various connectivity patterns. IBM's quantum systems employ error correction techniques and are accessible through cloud-based services.

Adiabatic Quantum Annealing (Algorithm \ref{alg:Annealing}) leverages principles of quantum superposition, tunneling, and coherence to find the ground state of a problem Hamiltonian. It starts with the system in the simple ground state of \( H_{\text{start}} \) and transitions to the target \( H_{\text{problem}} \) via a time-dependent Hamiltonian \( H(t) = A(t)H_{\text{start}} + B(t)H_{\text{problem}} \). If the transition is slow enough, the system remains in the ground state, achieving the optimal solution for the QUBO.
\begin{algorithm}[H]
\caption{Adiabatic Quantum Annealing for QUBO Problems}
\label{alg:Annealing}
\begin{algorithmic}[1]
\STATE $H_d \gets \text{Construct Driver Ising Hamiltonian}$
\STATE $|\psi(0)\rangle \gets \text{Prepare Ground State of} (H_d)$
\STATE $t \gets 0$
    
\WHILE{$t < T$}
    \STATE $s \gets t / T$
    \STATE $H(s) \gets (1-s)H_d + sH_p$
    \STATE $|\psi(t+\Delta t)\rangle \gets \text{Evolve State} (|\psi(t)\rangle, H(s), \Delta t)$
    \STATE $t \gets t + \Delta t$
\ENDWHILE
\STATE $|\psi_{\text{final}}\rangle \gets |\psi(T)\rangle$
\STATE $\text{solution} \gets \text{Measure State} (|\psi_{\text{final}}\rangle)$
\STATE \textbf{return} $\text{solution}$
\end{algorithmic}
\end{algorithm}
This study employs the quantum annealer developed by D-Wave Systems and the quantum-inspired solver from InfinityQ. The D-Wave Advantage processor, with over 5000 non-coherent qubits interconnected via the Pegasus topology, effectively addresses problems involving up to 180 fully-connected qubits. To handle larger-scale challenges, the D-Wave Hybrid algorithm combines classical and quantum approaches, breaking down complex problems into manageable subcomponents for efficient processing. In parallel, InfinityQ’s processor utilizes an Ising Machine approach, leveraging quantum-inspired techniques to solve large combinatorial optimization problems efficiently, thereby offering a robust alternative for tackling complex scenarios.

\section{Experiments}\label{sec:emp}

In these benchmarking experiments, we conducted numerical tests on both quantum and digital computing platforms to evaluate their efficiency and effectiveness in solving the portfolio optimization problems. We utilized a range of computational tools including CPU-based Binary Quadratic Programming (BQP) solvers by Gurobi with precise bounds, GPU-based heuristic QUBO solver qoqo/ABS2, and quantum annealing via D-Wave and InfinityQ, as well as quantum circuits through IBM Qiskit. The complexity of the optimization problem is influenced by several parameters: the number of assets \(n\), the number of time steps \(T\), the limit on asset selection \(B\), the number of blocks each asset can be chosen \(k\), the capital bound \(C\), and the risk-averse parameter \(q\). Among these, the parameters \(n\) and \(T\) dictate the difficulty of the problem. We examined the interplay and determined two combinations that pose greater challenges for both quantum and digital algorithms. Specifically, we designed one relatively smaller scale problem with \(n=200\) stocks and \(T=10\) that both digital and quantum algorithms with current hardware can solve, and one challenging large-scale problem with \(n=499\) and \(T=15\) that only digital can solve, establishing a benchmark for dynamic trading strategies in multi-period portfolio optimization under realistic settings in a QUBO framework:
\begin{itemize}
    \item[Exp1] \textbf{n200\_T10\_k3\_B60\_C10}: This scenario involves selecting from the top \(n=200\) market capitalization stocks in the S\&P 500 over two weeks (\(T=10\) trading days). Each asset can be selected up to \(k=3\) times in each period, with a maximum of \(B=60\) assets chosen and capital limited to \(C=10\) units. This period spans from October 2 to October 13, 2023, during which the S\&P 500 Market Index experienced a 0.47\% increase in the first week and a 0.18\% decrease in the second week. There are 12,100 variables in this QUBO setup.
    \item[Exp2] \textbf{n499\_T15\_k3\_B60\_C10}: This experiment extends to selecting from all \(n=499\) stocks of the S\&P 500 over three weeks (\(T=15\) trading days), with identical settings as the first experiment. The timeframe from April 29 to May 17, 2024, saw the S\&P 500 Market Index rising by 0.23\% in the first week, 0.81\% in the second week, and 1.57\% in the final week. This setup involves 45,060 variables in the QUBO framework.
\end{itemize}
The portfolio optimization problem incorporates short selling, transaction costs, and path-dependent cash constraints. Table~\ref{tab:experiment_parameters} outlines the parameters utilized in our experiments, chosen based on relevant literature. It includes a transaction cost of 0.1\% of the traded value, a daily risk-free rate of 0.01\%, and a short sale loan rate of 0.0025\% \citep{gaivoronski2005optimal, rietz1988equity, d2002market, tepla2000optimal}.

\begin{table}[h]
\centering
\begin{tabular}{l|l}
\hline\hline
\textbf{Parameter} & \textbf{Value} \\
\hline
Daily risk-free rate, \( \rho_c \) & 0.01\% (annual rate = 2.55\%) \\
Transaction cost rate, \( \delta \) & 0.1\% \\
Daily loan rate on short positions, \( \rho_s \) & 0.0025\% (annual rate = 0.92\%) \\
Available capital \(C\) & 10 units (\$1,000,000) \\
Maximum \# assets at each time step \(B\) & 60 \\
Maximum \# blocks each asset \(k\) & 3 \\
\hline\hline
\end{tabular}
\caption{Parameters for Portfolio Optimization}
\label{tab:experiment_parameters}
\end{table}
The risk-averse parameter \(q\) significantly influences the problem's solvability. For two extremes, the multi-objective problem reduces to single objective and is easier to solve. Specifically, when \(q = 0\), the investor is completely profit-oriented, ignoring risk, while as \(q \to \infty\), the investor becomes extremely risk-averse, preferring to hold cash as the risk-free asset. On the contrary, the optimization problem becomes more difficult to solve, where a balanced approach considering different levels of risk preference is adopted. We analyze a range of eight values for \(q\), from \(q = 0\) (risk-seeking) to \(q = 0.01\) (near risk-free) and show the Pareto optimal frontier, aiming to provide a clear understanding of the current capabilities of modern computing technologies in addressing real-world financial scenarios.

\subsection{Evaluation Criteria}
To assess the performance of both quantum and digital computing platforms in solving portfolio optimization problems, we employ two primary types of evaluation criteria: computational and economic perspectives.

\begin{itemize}
    \item \textbf{Computational Performance}: This category assesses the technical quality of the solutions:
    \begin{itemize}
        \item \textbf{Objective Value}: Ideally, lower values are better as they indicate a closer alignment with the defined financial goals.
        \item \textbf{Optimality Gap}: This measures the difference between the solver's best-known solution and the theoretical optimum (e.g., lower bound (LB) of BQP if integer constraints were relaxed), offering insights into possible improvements. Note that reaching the LB might not always result in feasible solutions under integer constraints. The gap is calculated as:
\begin{equation*}
\text{Gap} = \frac{\lvert \text{Solver's Obj} - \text{BQP's LB} \rvert}{\lvert \text{Solver's Obj} \rvert}
\end{equation*}
        \item \textbf{Time-to-Solution (TTS)}: This metric records the time taken by the algorithm to first find the best feasible solution, highlighting computational efficiency. This measure is crucial when objective values are comparable, as it provides a direct comparison of the speed and efficacy of quantum algorithms versus traditional classical algorithms.
    \end{itemize}
    \item \textbf{Economic Performance}: These criteria reflect the economic viability and practical implementation of the trading strategies:
    \begin{itemize}
        \item \textbf{Volatility}: Lower volatility of portfolio is preferred, as it indicates a more stable and less risky investment, aligning with the goal of risk minimization.
        \item \textbf{Profits}: Higher profits signify a successful investment strategy, marking the solver's ability to maximize financial returns.
        \item \textbf{Trading Strategy Trajectory}: Analyzes the sequence of trades executed during the investment period. In scenarios like \(q = 0.01\), where risk aversion is extreme, holding cash should be optimal, serving as a validity check for the solver's strategic output. In balanced scenarios, a sparse portfolio is often economically and practically preferable.
        \item \textbf{Sharpe Ratio}: This ratio assesses investment performance by adjusting for risk, facilitating quick and direct comparisons across investments. A higher Sharpe ratio denotes a more attractive risk-adjusted return.
        \item \textbf{Transaction Costs}: Minimizing these costs is essential as they directly reduce net profits. Effective trading strategies should aim for high returns while managing and minimizing associated transactional expenses.
    \end{itemize}
\end{itemize}
These metrics collectively provide a holistic view of the computational effectiveness and economic impact of the optimization solutions, enabling a comprehensive evaluation of solvers on their practical viability and financial soundness.

\subsection{Findings and Implications}
The experiments investigate the computational performance of digital computing -- BQP solver Gurobi -- and quantum-ready QUBO solver qoqo/ABS2 as well as quantum computing platforms including D-Wave, IBM Qiskit, and InfinityQ. While quantum computing is rapidly evolving and faces challenges, particularly in hardware robustness and noise correction, the landscape is promising. The inherent compatibility of QUBO with quantum algorithms ensures that as more powerful and error-corrected quantum computers become available, the transition will be seamless. Therefore, in the comparative analysis, we focus on benchmarking quantum-ready algorithms (QUBO) against their classical counterparts (BQP), anticipating significant advancements as quantum computing matures.

\subsubsection{Experiment 1: n200\_T10\_k3\_B60\_C10}
The experiment aims to optimize a medium-sized portfolio allocation among 200 assets and cash over 10 trading days. BQP problems were solved using Gurobi v11.0 with a 10,000-second time limit, while QUBO problems were addressed using ABS2 within a 3,600-second limit. 
\begin{table}[ht]
\centering
\renewcommand{\arraystretch}{1.1}
\resizebox{\textwidth}{!}{
\begin{tabular}{l|r|r|r|c|r|r|c}
\hline\hline
\hspace{3em} $q$ & \textbf{LB} & \textbf{UB} & \textbf{Gap} & \textbf{Gurobi TTS [s]} & \textbf{QUBO Obj.} & \textbf{Gap*} & \textbf{ABS2 TTS [s]} \\ \hline
0 & -1,855,197 & -1,855,197 & 0.000\% & 0.22 & -1,855,197 & 0.000\% & 102 \\ \hline
0.000001 & -1,839,242 & -1,839,242 & 0.000\% & 34 & -1,839,012 & 0.010\% & 210 \\ \hline
0.00001 & -1,706,529 & -1,697,030 & 0.560\% & 125 & -1,701,223 & 0.312\% & 105 \\ \hline
0.00005 & -1,289,791 & -1,265,722 & 1.891\% & 3,212 & -1,269,382 & 1.608\% & 1,191 \\ \hline
0.0001 & -988,003 & -950,586 & 3.936\% & 9,591 & -953,432 & 3.626\% & 1,778 \\ \hline
0.0005 & -356,010 & -296,540 & 20.055\% & 10,000 & -272,033 & 30.870\% & 3,562 \\ \hline
0.001 & -231,724 & -168,674 & 37.380\% & 10,000 & -104,856 & 120.992\% & 3,552 \\ \hline
0.01 & -108,821 & -1,000 & 10,782\% & 33 & -1,000 & 10,782\% & 53 \\ \hline\hline
\end{tabular}
}
\caption{BQP (classical) vs QUBO (quantum-ready) in portfolio optimization of n200\_T10\_k3\_B60\_C10}
\label{tab:a200_t10}
\end{table}
As shown in Table~\ref{tab:a200_t10}, the computational performance varies with the risk aversion parameter \(q\). Notably, as \(q\) approaches extremes (either focusing solely on profit or risk), the problem becomes significantly easier to solve, with solvers quickly finding solutions. Particularly, for extreme risk aversion (\(q=0.01\)), both solvers converge to a solution in under a minute, albeit with a considerable gap indicative of the distance from the theoretical optimum under relaxed constraints. 

\begin{figure}[ht]
    \centering
    \includegraphics[width=0.8\textwidth]{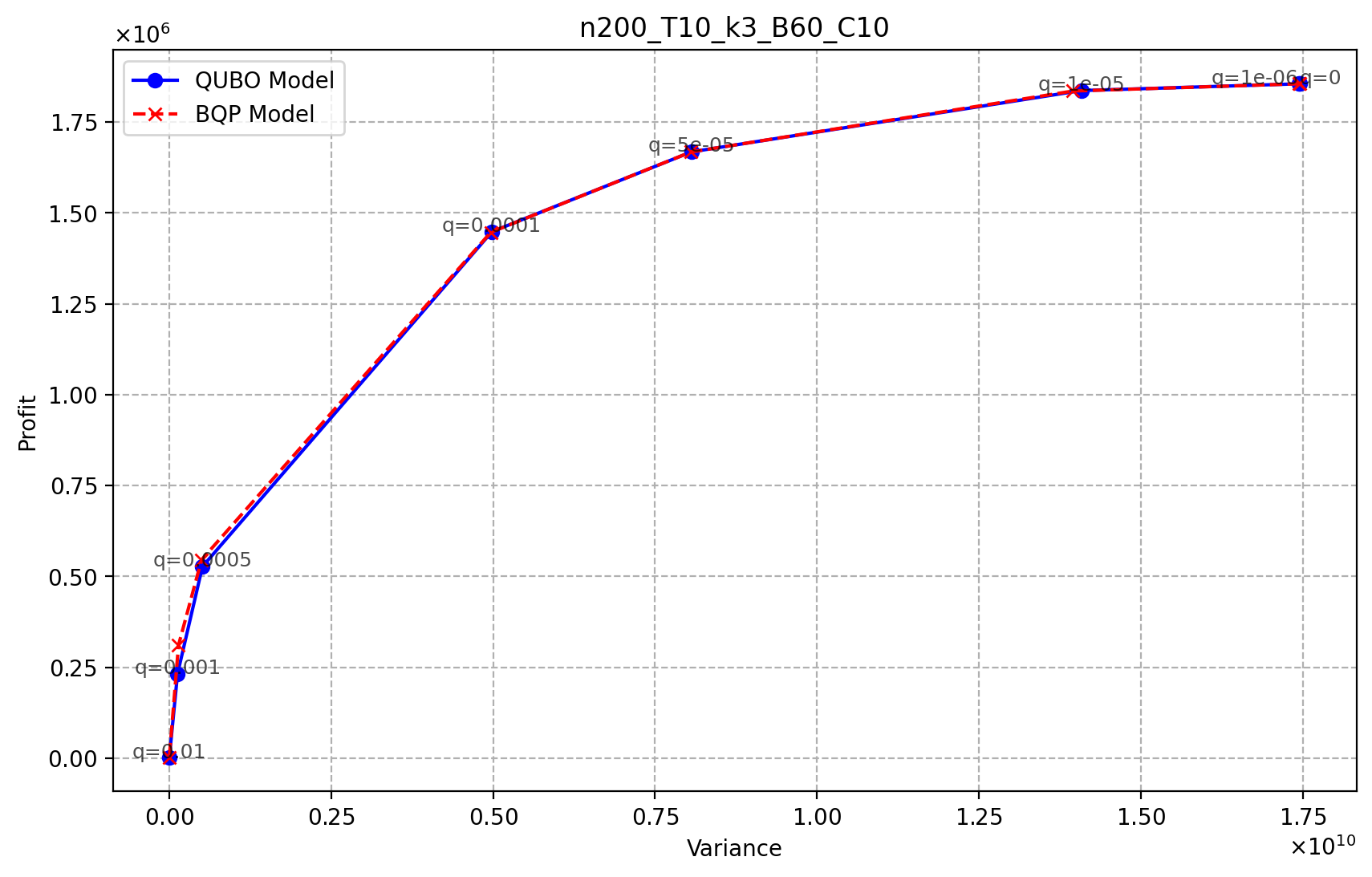}
    \caption{Pareto optimal frontier for n200\_T10\_k3\_B60\_C10}
    \label{fig:profit_var_200}
\end{figure}

\begin{figure}[h!]
    \centering
    \begin{subfigure}[b]{0.47\textwidth}
        \includegraphics[width=\textwidth]{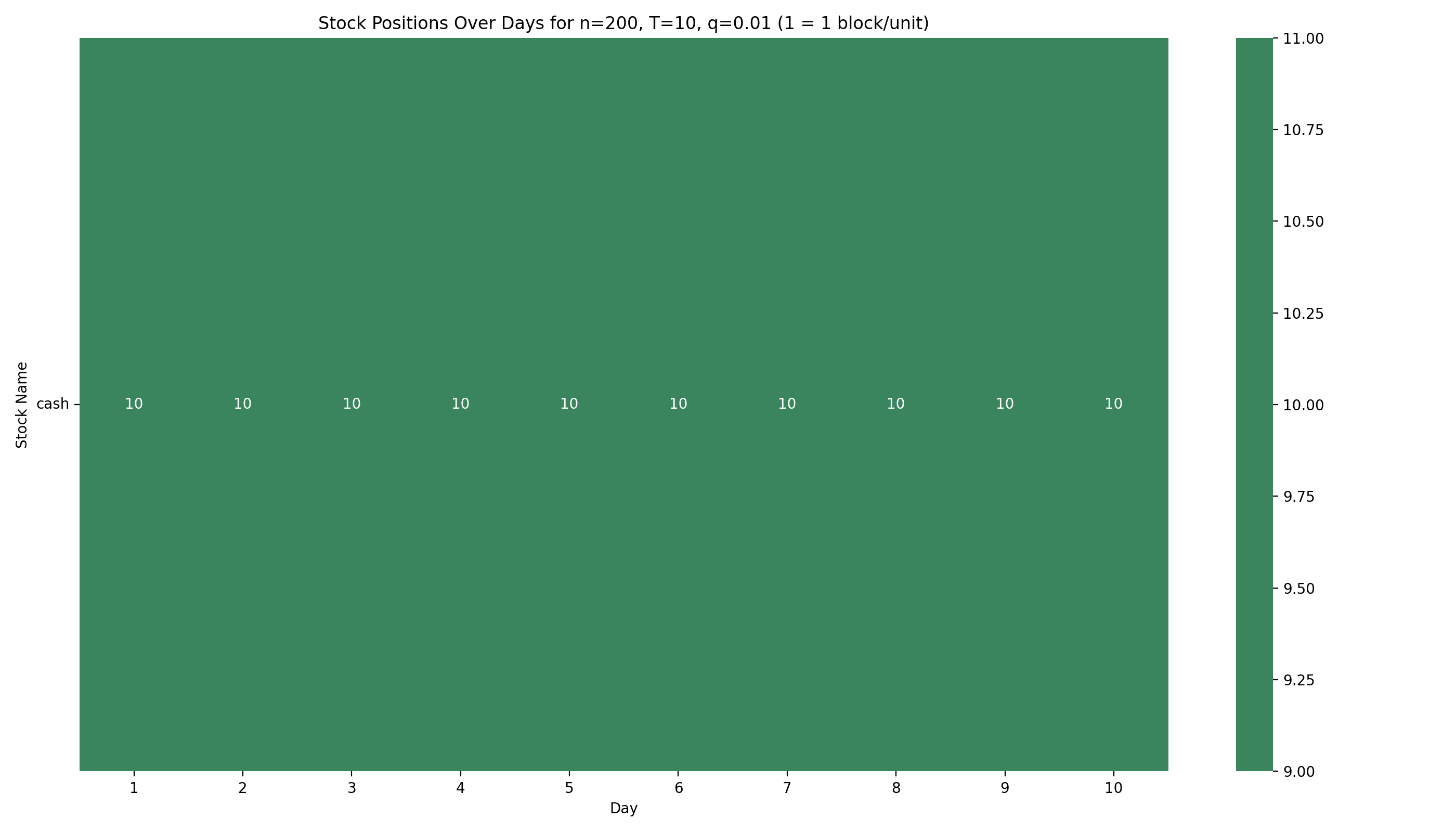}
        \caption{BQP model strategy (Gurobi)}
    \end{subfigure}
    \hfill
    \begin{subfigure}[b]{0.47\textwidth}
        \includegraphics[width=\textwidth]{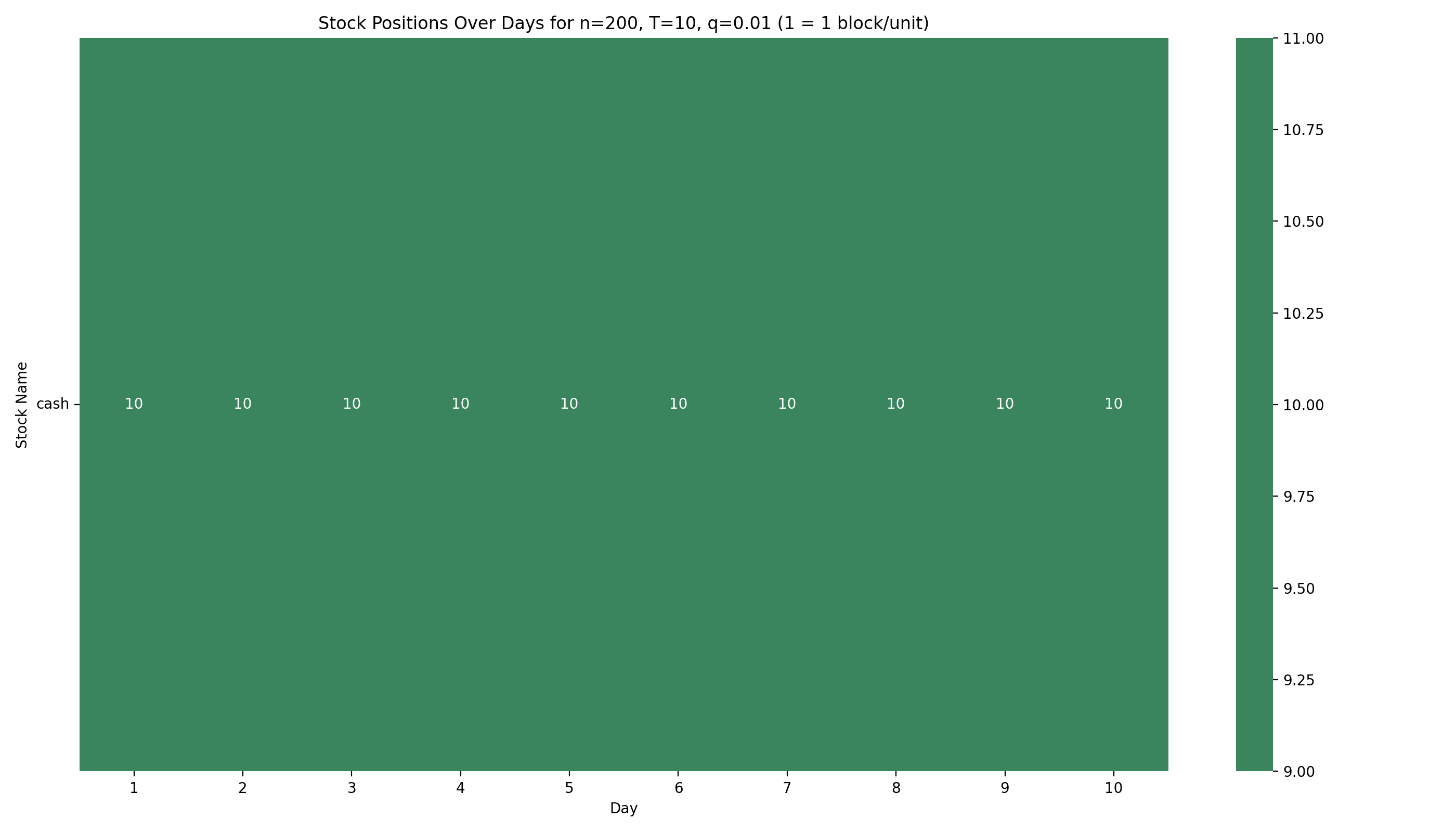}
        \caption{QUBO model strategy (ABS2)}
    \end{subfigure}
    \caption{Stock positions for n200\_T10\_k3\_B60\_C10 with $q = 0.01$}
    \label{fig:n200_T10}
\end{figure}

\begin{figure}[h!]
    \centering
    \begin{subfigure}[b]{0.47\textwidth}
        \includegraphics[width=\textwidth]{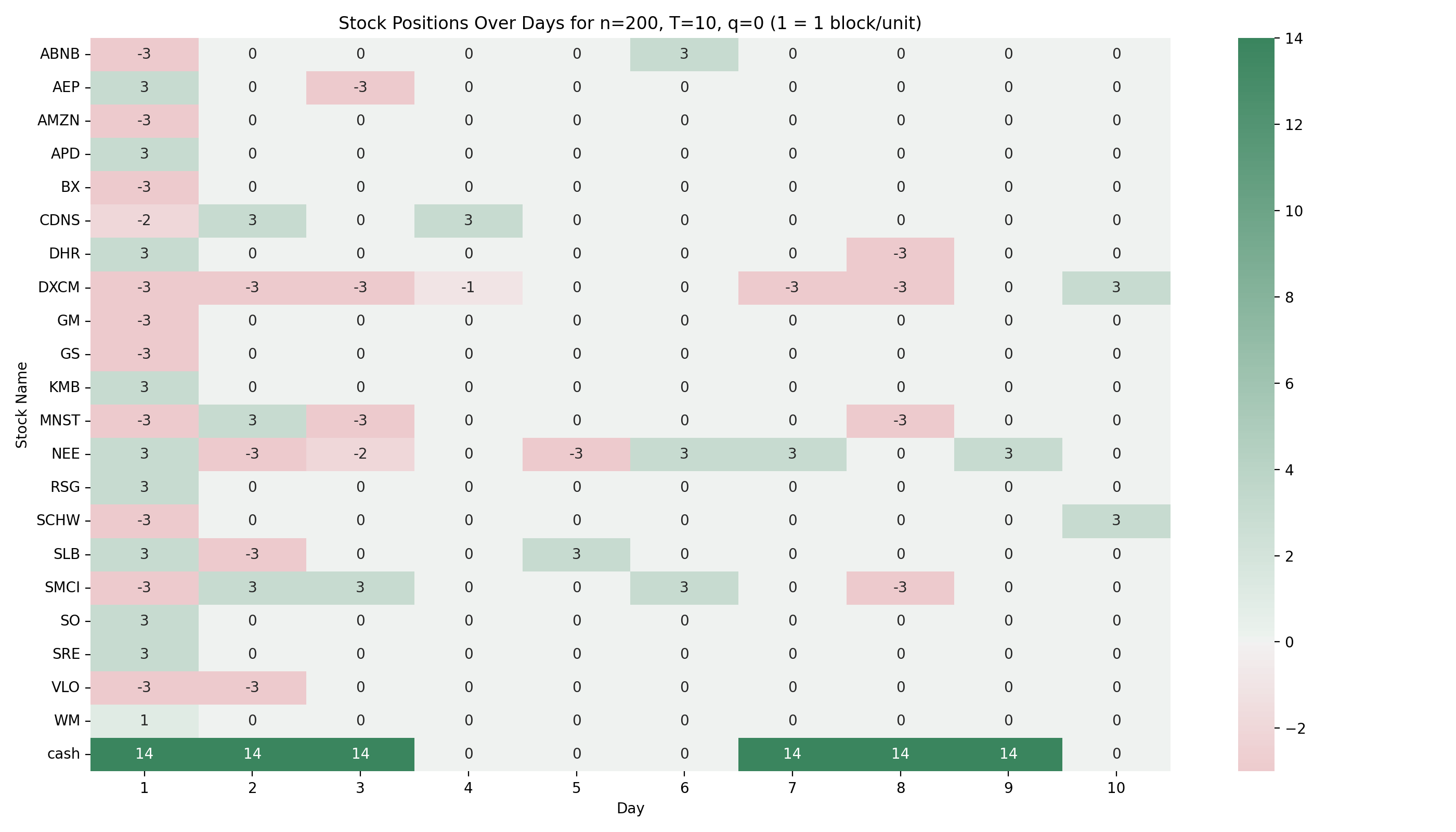}
        \caption{BQP model strategy (Gurobi)}
        \label{fig:sub_bqp}
    \end{subfigure}
    \hfill
    \begin{subfigure}[b]{0.47\textwidth}
        \includegraphics[width=\textwidth]{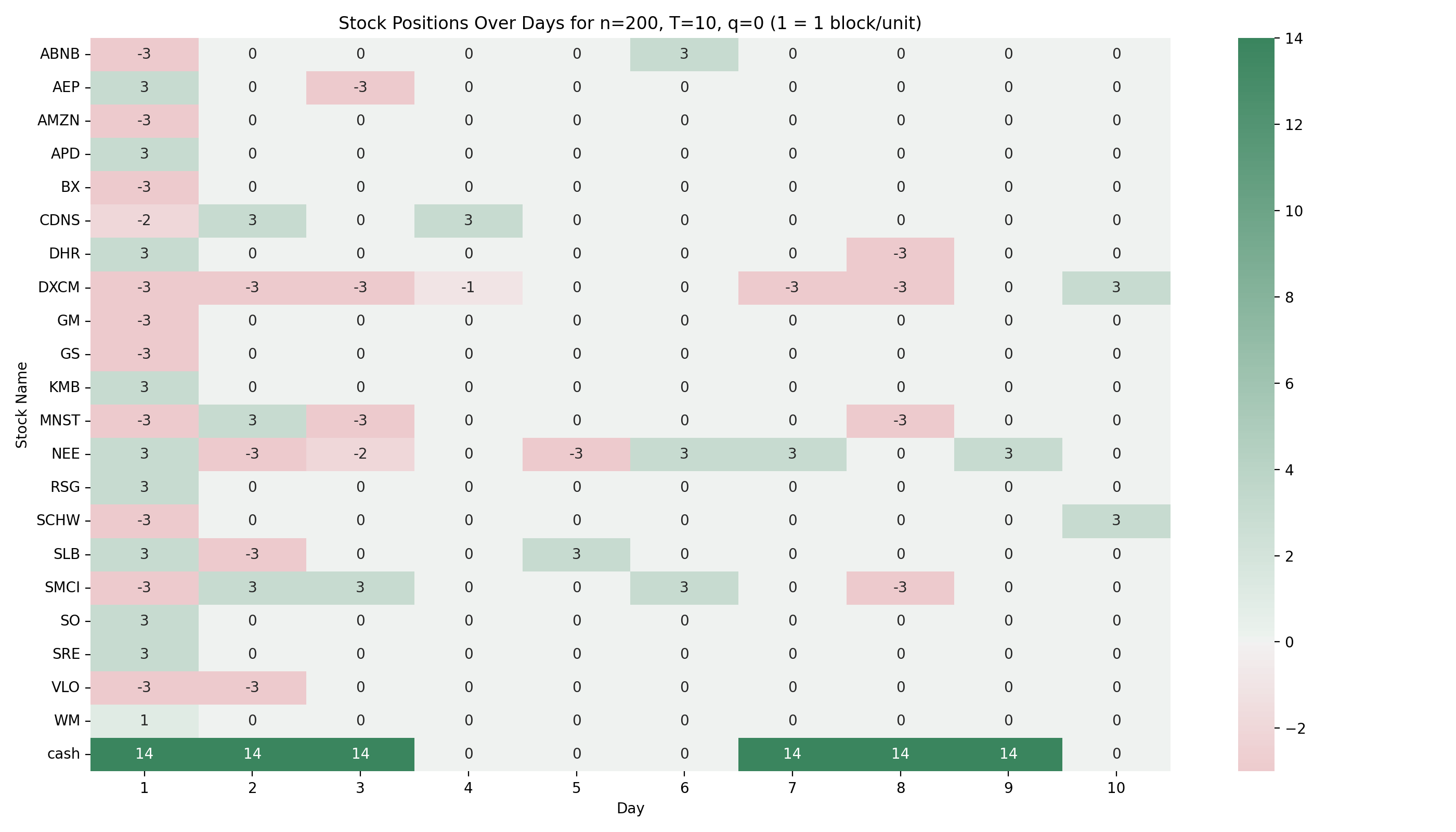}
        \caption{QUBO model strategy (ABS2)}
        \label{fig:sub_qoqo}
    \end{subfigure}
    \caption{Stock positions for n200\_T10\_k3\_B60\_C10 with \(q = 0\)}
    \label{fig:position_200_q0}
\end{figure}

Figure~\ref{fig:profit_var_200} shows the Pareto optimal frontier, aligning with the computed results and showing similar profit-risk dependence against various values of \(q\).
When \(q\) is large enough, e.g., at \(q=0.01\), it effectively minimizes risk to nearly zero. Both methods can solve the problem within one minute, but exhibit an extremely large gap of 10,782\%. However, the gap is merely an indicator reflecting the distance to the theoretical optimum under relaxed constraints and does not necessarily imply that the solutions are suboptimal. In this extreme risk-averse scenario, the best strategy should be holding cash over time. Both the BQP and QUBO strategies consistently recommend holding cash up to the maximum of \(C=10\) blocks throughout the period, as illustrated in Figure~\ref{fig:n200_T10}. This consistency validates that the solutions are globally optimal.
For \(q = 0\) up to \(q = 0.00001\), risk considerations become irrelevant, and the focus shifts entirely towards maximizing profit. Both Gurobi and ABS2 efficiently find solutions within 2 minutes, closely matching or achieving the theoretical optimal benchmark under relaxed constraints (gap close to or equal to 0\%). Up to \(q = 0.0001\), the problem remains solvable within 2.5 hours with a reasonable gap of less than 5\%. Figure~\ref{fig:sub_bqp} and~\ref{fig:sub_qoqo} illustrate the stock holdings and cash positions of the two models. Due to space constraints, the figures only display the stocks selected in the first time period (\(t=1\)) and their corresponding positions over time up to \(T\). For both models, the dynamic positions are identical.

\begin{figure}[ht]
    \centering
    \begin{subfigure}[b]{0.47\textwidth}
        \centering
        \includegraphics[width=\textwidth]{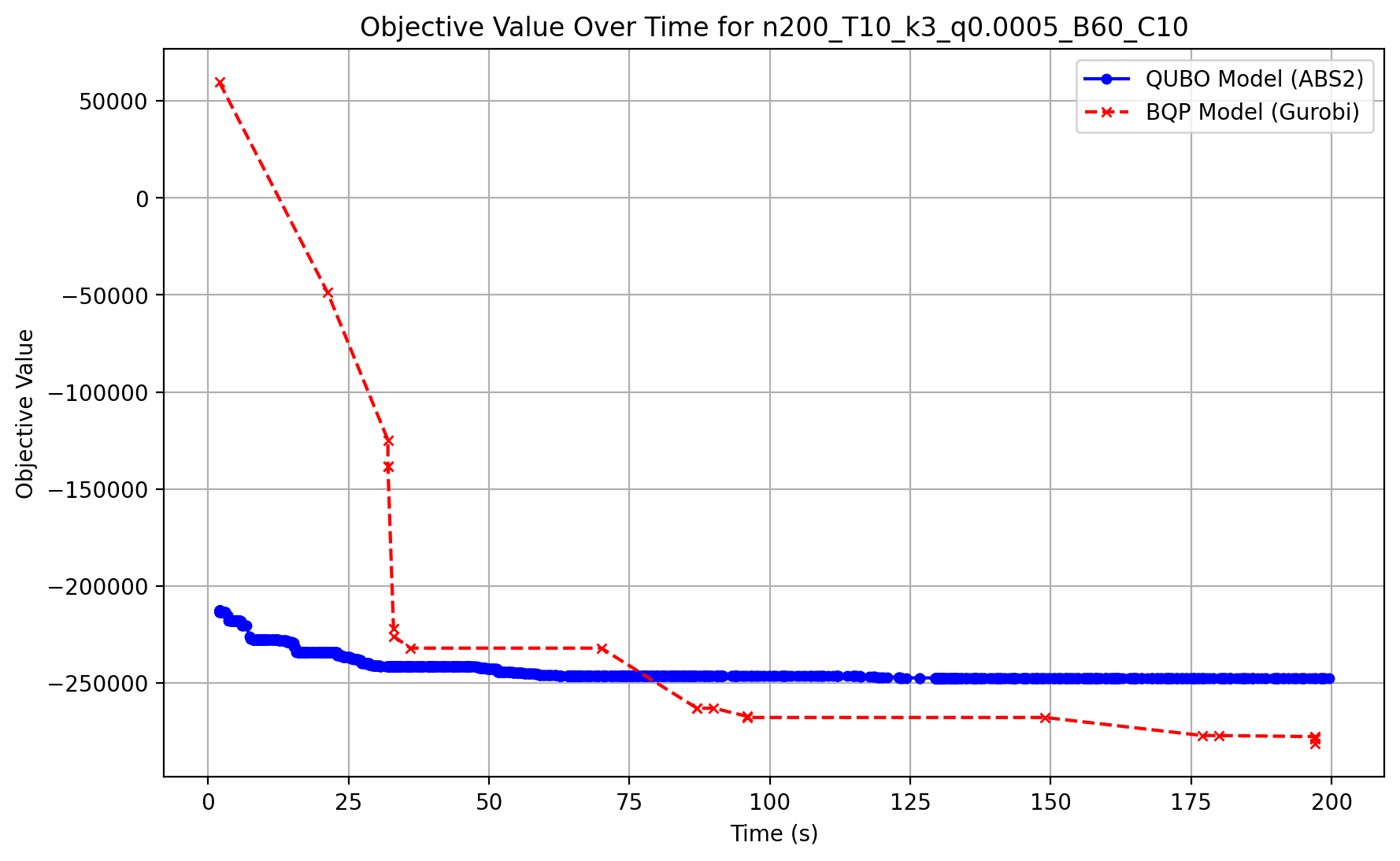}
        \caption{Objective Change for \( q = 0.0005 \) in a200.}
        \label{fig:a200_t10_q0.0005}
    \end{subfigure}
    \hfill
    \begin{subfigure}[b]{0.47\textwidth}
        \centering
        \includegraphics[width=\textwidth]{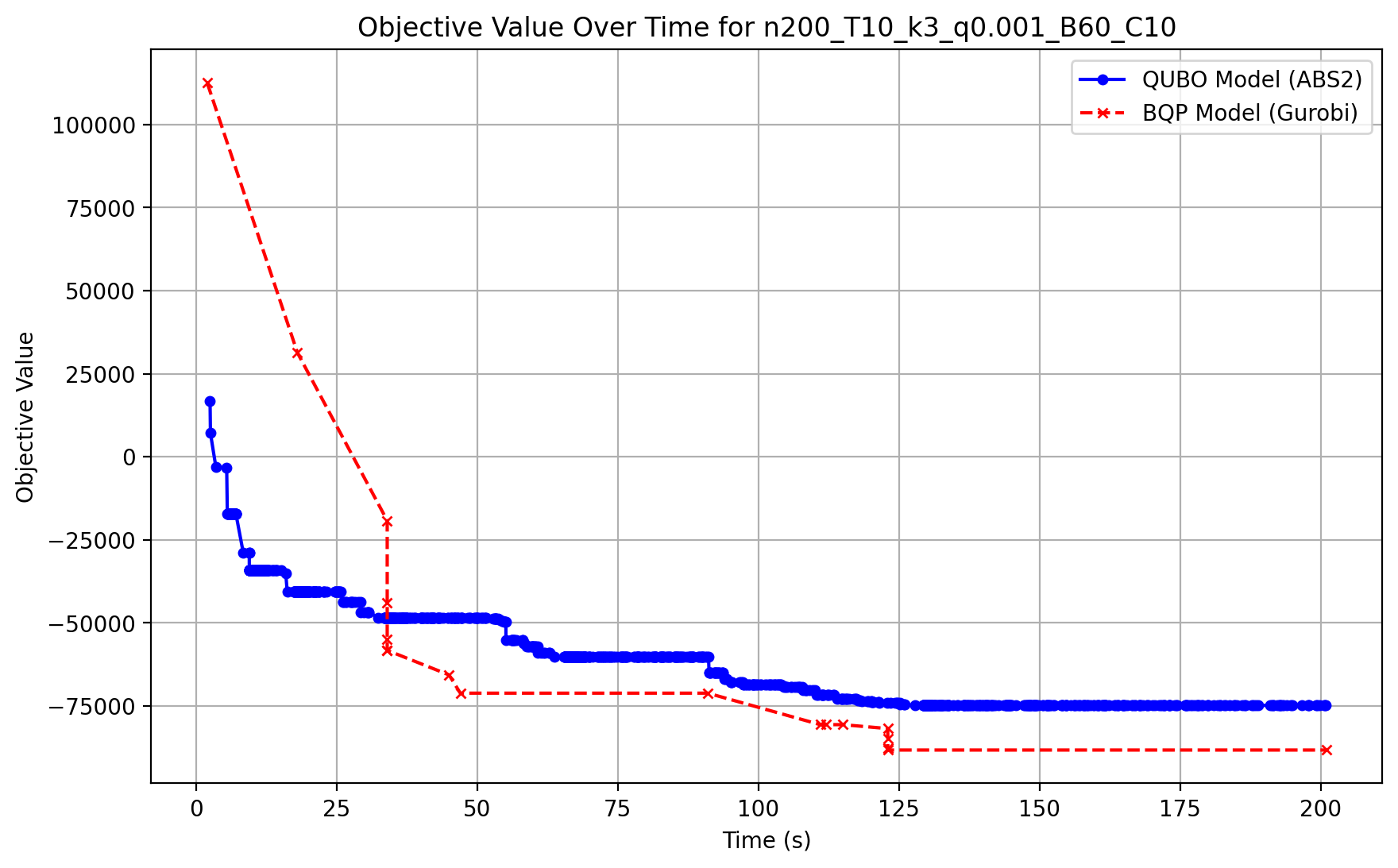}
        \caption{Objective Change for \( q = 0.001 \) in a200.}
        \label{fig:a200_t10_q0.001}
    \end{subfigure}
    \caption{Convergence of objective values for different risk aversion parameters.}
\end{figure}

While Gurobi generally finds solutions faster than qoqo/ABS2 for simpler problems, in scenarios where \(q\) is set to manage both profit and risk, the optimization problem becomes challenging, demanding more computational resources. For \(q\) values such as 0.0005 and 0.001, both solvers nearly exhaust the computational limits set. Even though Gurobi generally provides a slightly better solution given longer computational time limit, it requires substantially more time to reach to a good solution with optimal \(\pm\epsilon\) compared to qoqo ABS2, as detailed in the convergence plots, see Figures~\ref{fig:a200_t10_q0.0005} and~\ref{fig:a200_t10_q0.001}. Specifically, ABS2 converges towards a good solution significantly faster than Gurobi, e.g., within 71 seconds for \( q = 0.0005 \) and within 35 seconds for \( q = 0.001 \). Both solvers cannot significantly improve upon them further until the computational time limit is reached.
When we further vary the computational time limits to 120, 600, 1800, and 3600 seconds, we observe that at a very short time limit of 120 seconds, the objective values of QUBO and BQP are very close for different \( q \) values. At 600, 1800, and 3600 seconds, except for the most challenging problem at \( q = 0.001 \), where the objective value of QUBO optimized by ABS2 is noticeably higher than that of BQP optimized by Gurobi, the differences between the objective values of QUBO and BQP for \( q = 0.001 \) decrease as the time limit increases.

\begin{figure}[h!]
    \centering
    \begin{subfigure}[b]{0.47\textwidth}
        \includegraphics[width=\textwidth]{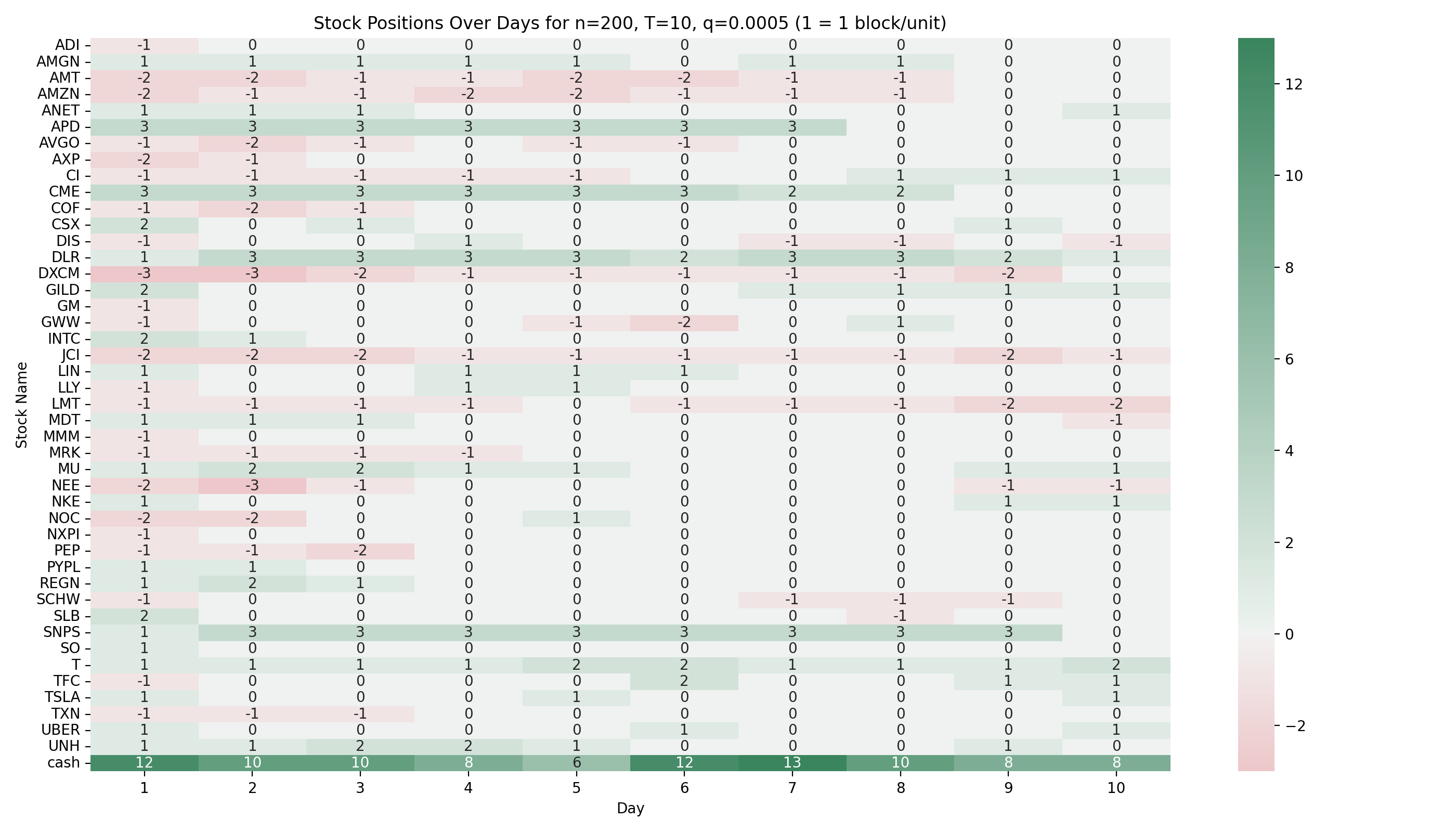}
        \caption{BQP model strategy (Gurobi)}
        \label{fig:bqp_200_q0005}
    \end{subfigure}
    \hfill
    \begin{subfigure}[b]{0.47\textwidth}
        \includegraphics[width=\textwidth]{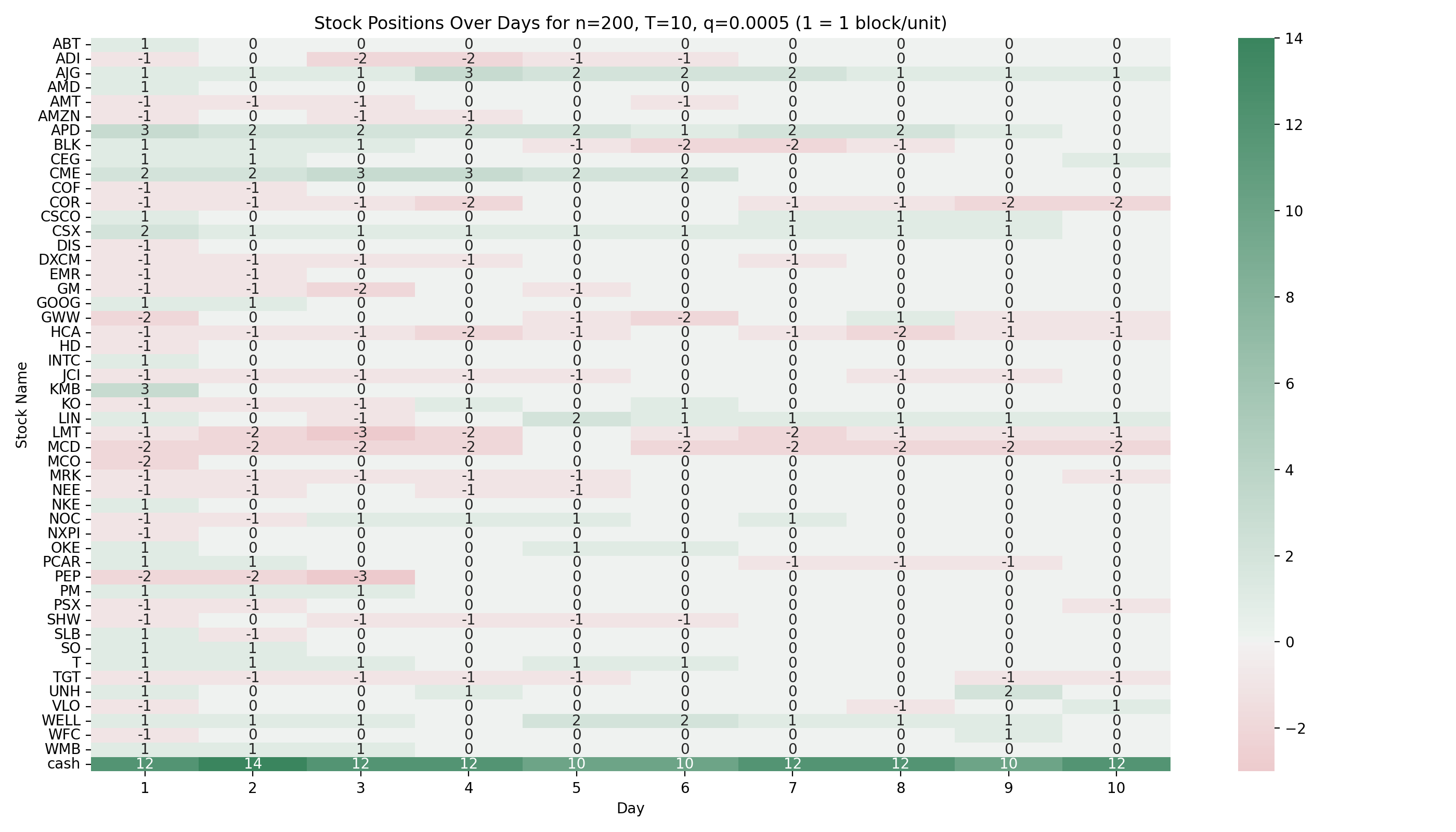}
        \caption{QUBO model strategy (ABS2)}
        \label{fig:qubo_200_q0005}
    \end{subfigure}
    \caption{Stock positions for n200\_T10\_k3\_B60\_C10 with $q = 0.0005$}
\end{figure}

Figures \ref{fig:bqp_200_q0005} and \ref{fig:qubo_200_q0005} display the trading trajectories at \( q = 0.0005 \), where both risk and profit are considered. The two models provide different strategies. On the first day, the BQP model selected 44 stocks, while the QUBO model selected 50 stocks. Nevertheless, their positions on the same stocks over the period were quite similar. Specifically, both models held long positions in APD from day 1 to day 7, long positions in CME from day 1 to day 8, short positions in LMT on all days except the fifth, and short positions in PEP for the first three days.

\subsubsection{Experiment 2: n499\_T15\_k3\_B60\_C10}
In this larger-scale experiment, we address the portfolio optimization problem involving the allocation of capital among 499 stocks and cash over a period of three weeks while keeping other parameters the same. Given the current limitations of quantum computing hardware, implementing quantum solutions for such extensive and complex scenarios remains out of reach. However, we establish a robust benchmark using digital computing as a reference point. This benchmark is designed with the anticipation of rapid advancements in quantum computing technologies. Our goal is to provide a transparent and fair platform to evaluate genuine breakthroughs in quantum computing as they emerge, helping to determine when quantum solutions can realistically surpass the capabilities of traditional digital algorithms in handling large-scale financial optimization tasks.

\begin{figure}[ht]
    \centering
    \begin{subfigure}[b]{0.8\textwidth}
        \centering
        \includegraphics[width=\textwidth]{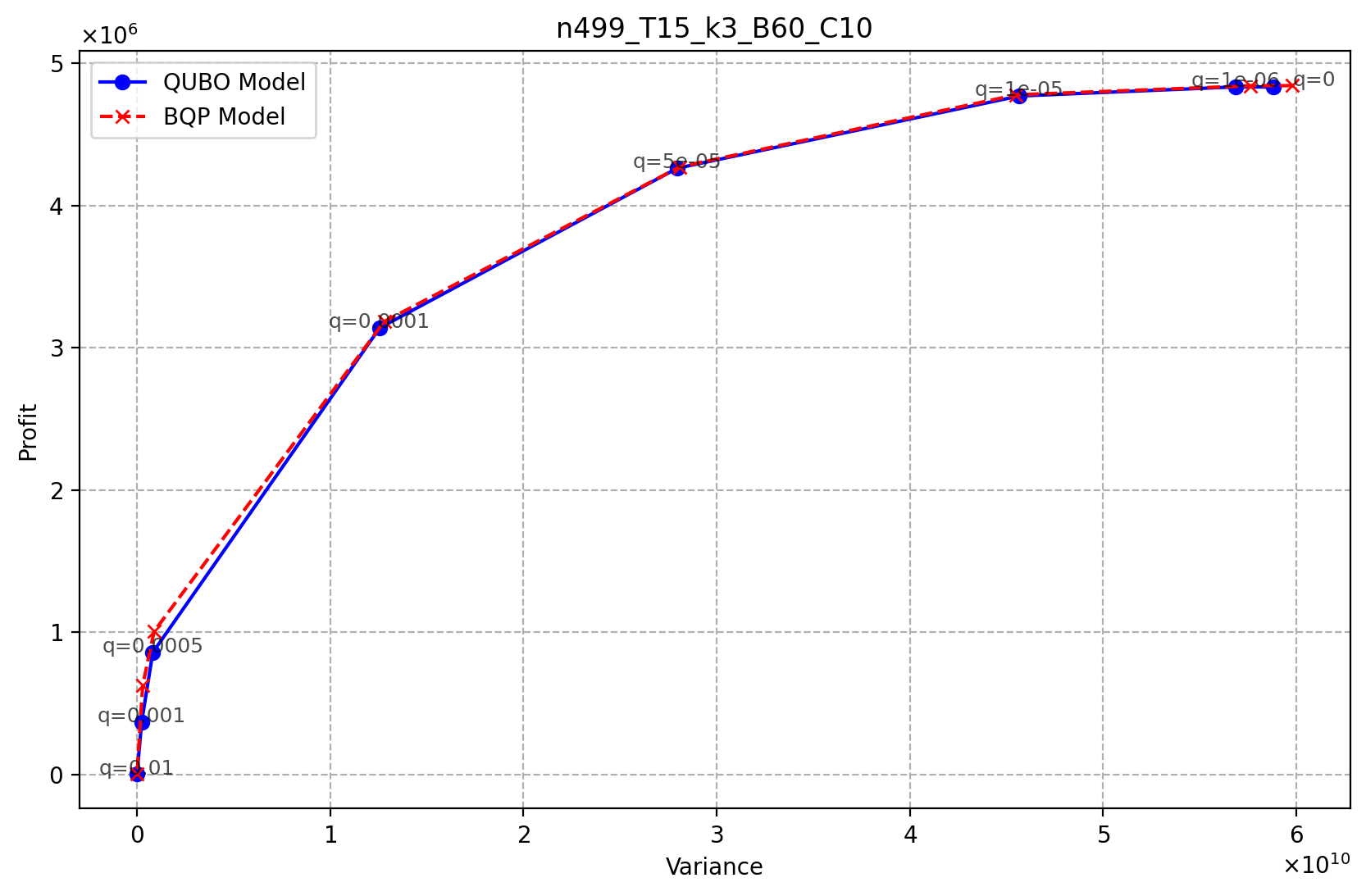}
        \caption{n499\_T15\_k3\_B60\_C10}
        \label{fig:profit_var_499}
    \end{subfigure}
    \caption{Profits and Variances at different \( q \) values.}
\end{figure}

For this experiment, BQP problems were addressed using Gurobi v11.0, with a time limit of 10,000 seconds, while QUBO problems were solved using the ABS2 solver, with a time limit of 3,600 seconds. The performance comparison between these digital computing approaches in portfolio optimization is presented in Table~\ref{tab:a499_t15}. A Pareto optimal frontier, displayed in Figure~\ref{fig:profit_var_499}, illustrates that the solvers provide comparable results under various scenarios.
\begin{table}[ht]
\centering
\renewcommand{\arraystretch}{1.1}
\resizebox{\textwidth}{!}{
\begin{tabular}{l|r|r|r|c|r|r|c}
\hline\hline
\hspace{3em} $q$     & \textbf{LB}      
 & \textbf{UB}      & \textbf{Gap}  & \textbf{Gurobi TTS [s]} & \textbf{QUBO Obj.} & \textbf{Gap*} & \textbf{ABS2 TTS [s]} \\ \hline
0        & -4,841,118 & -4,841,118 & 0.000\%  & 0.93           & -4,831,302  & 0.203\%   & 3,615        \\ \hline
0.000001 & -4,795,977 & -4,785,071 & 0.228\%  & 588             & -4,787,080  & 0.186\%   & 3,568        \\ \hline
0.00001  & -4,343,457 & -4,324,080 & 0.448\%  & 834            & -4,322,507  & 0.484\%   & 3,621        \\ \hline
0.00005  & -2,927,666 & -2,867,493 & 2.099\%  & 3,060          & -2,871,591  & 1.953\%   & 3,606      \\ \hline
0.0001   & -1,996,723   & -1,903,969   & 4.872\%  & 2,663          & -1,887,560    & 5.783\%   & 3,545      \\ \hline
0.0005   & -708,516   & -563,441   & 25.748\% & 4,843         & -451,312    & 56.990\%  & 3,579      \\ \hline
0.001    & -507,219   & -332,172   & 52.698\% & 5,243         & -132,613    & 282.481\% & 3,595      \\ \hline
0.01     & -316,320   & -1,500     & 20,988\% & 337             & -1,500      & 20,988\%  & 1,026         \\ \hline\hline
\end{tabular}
}
\caption{BQP vs QUBO in portfolio optimization of n499\_T15\_k3\_B60\_C10}
\label{tab:a499_t15}
\end{table}
As the risk aversion parameter \(q\) approaches extremities -- focusing solely on profit or risk, the problem is easier to solve. When \(q \rightarrow 0\), the small optimality gap verifies the effectiveness of the solutions; whereas at a high \(q\), e.g., \(q=0.01\), both solvers suggest holding cash, effectively minimizing risk to zero and thereby justifying the optimality of the solutions despite a large gap. Generally, Gurobi offers better efficiency with shorter TTS for simpler problems, while QUBO resolves more challenging scenarios at \(q=0.0005\) and \(q=0.001\) more effectively.

\begin{figure}[ht]
    \centering
    \begin{subfigure}[b]{0.47\textwidth}
        \centering
        \includegraphics[width=\textwidth]{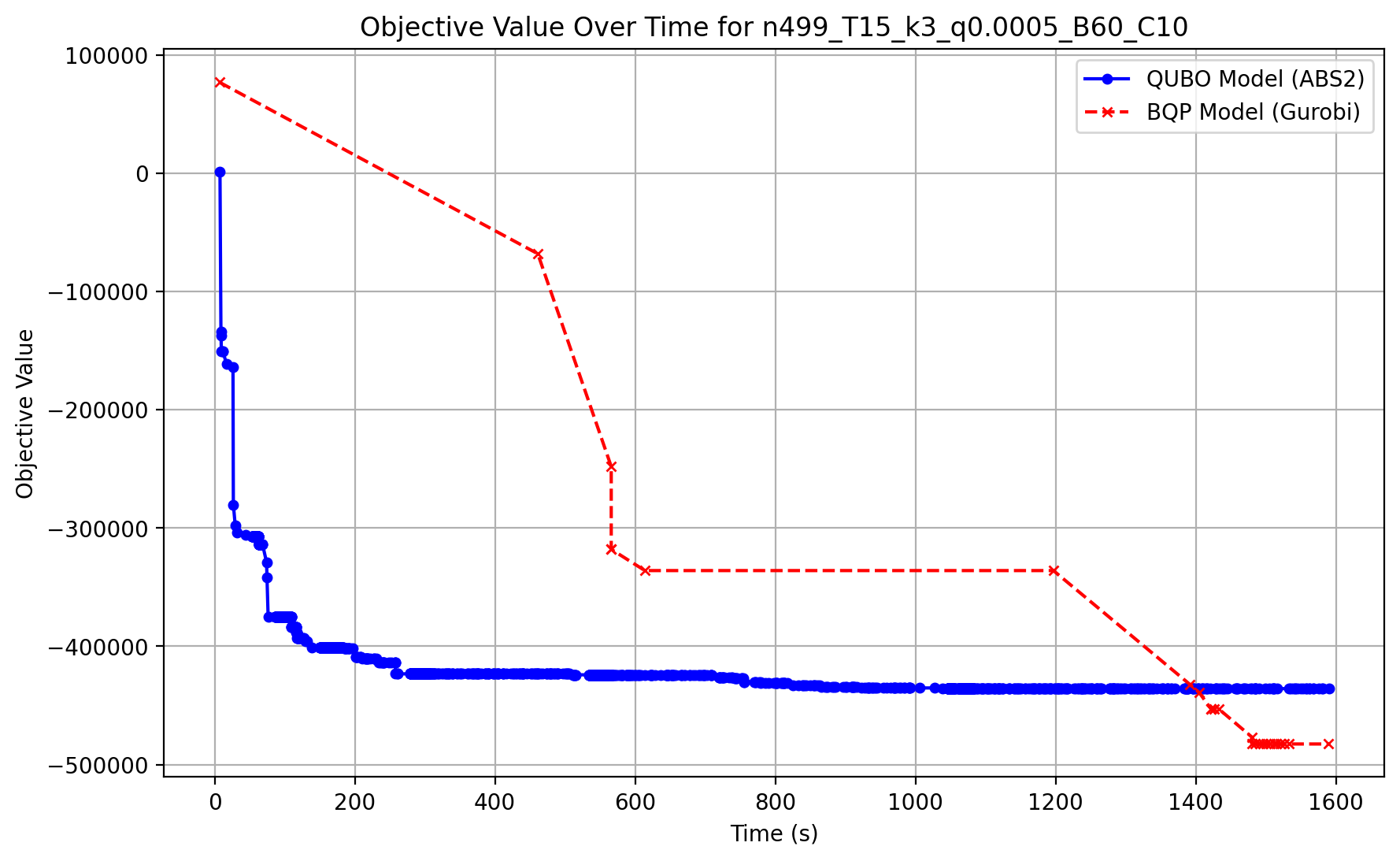}
        \caption{Objective Change for \( q = 0.0005 \) in a499.}
        \label{fig:a499_t15_q0.0005}
    \end{subfigure}
    \hfill
    \begin{subfigure}[b]{0.47\textwidth}
        \centering
        \includegraphics[width=\textwidth]{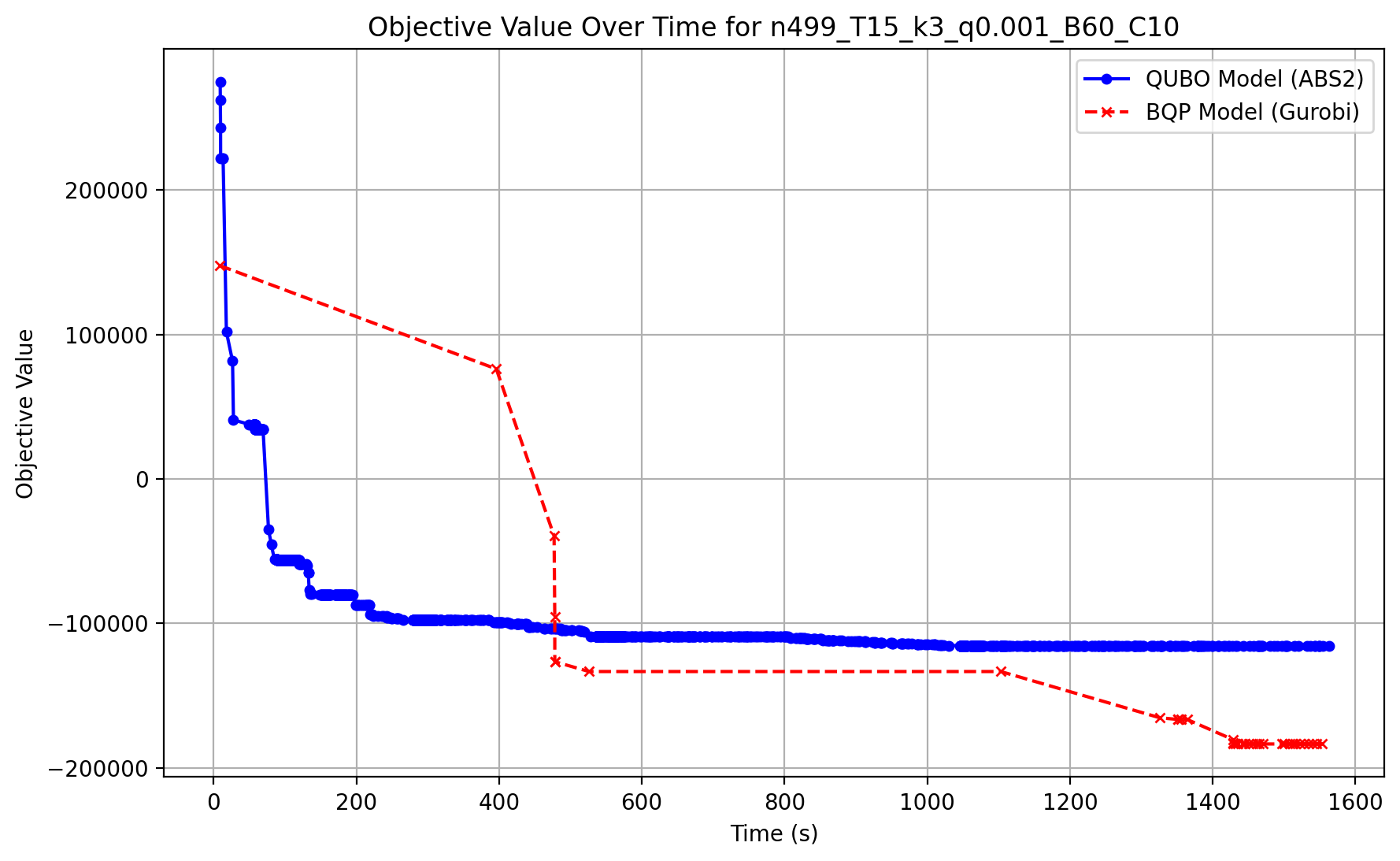}
        \caption{Objective Change for \( q = 0.001 \) in a499.}
        \label{fig:a499_t15_q0.001}
    \end{subfigure}
    \caption{}
\end{figure}

\begin{figure}[ht]
    \centering
    \begin{subfigure}[b]{0.47\textwidth}
        \centering
        \includegraphics[width=\textwidth]{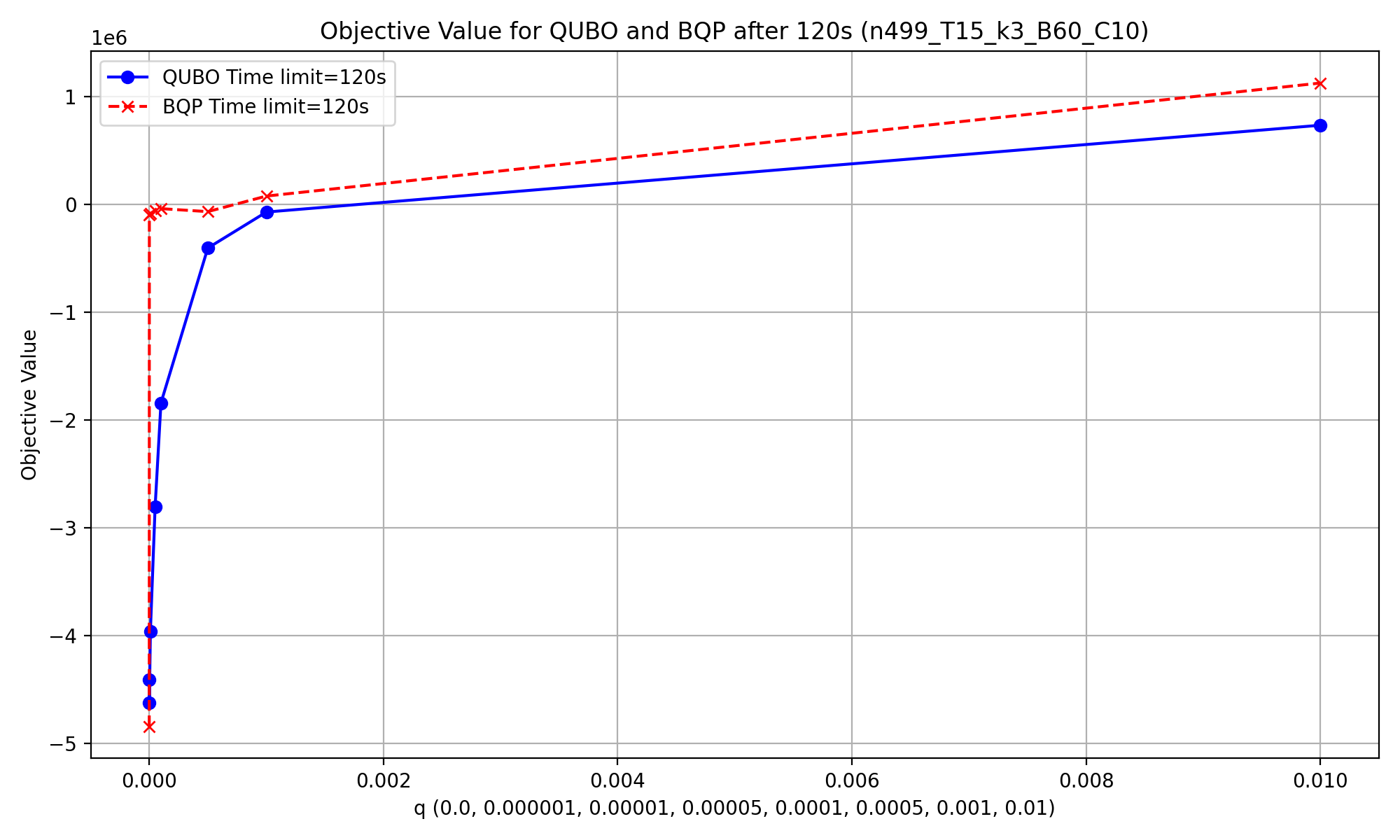}
        \caption{}
        \label{fig:comparison_a499_t15_120s}
    \end{subfigure}
    \hfill
    \begin{subfigure}[b]{0.47\textwidth}
        \centering
        \includegraphics[width=\textwidth]{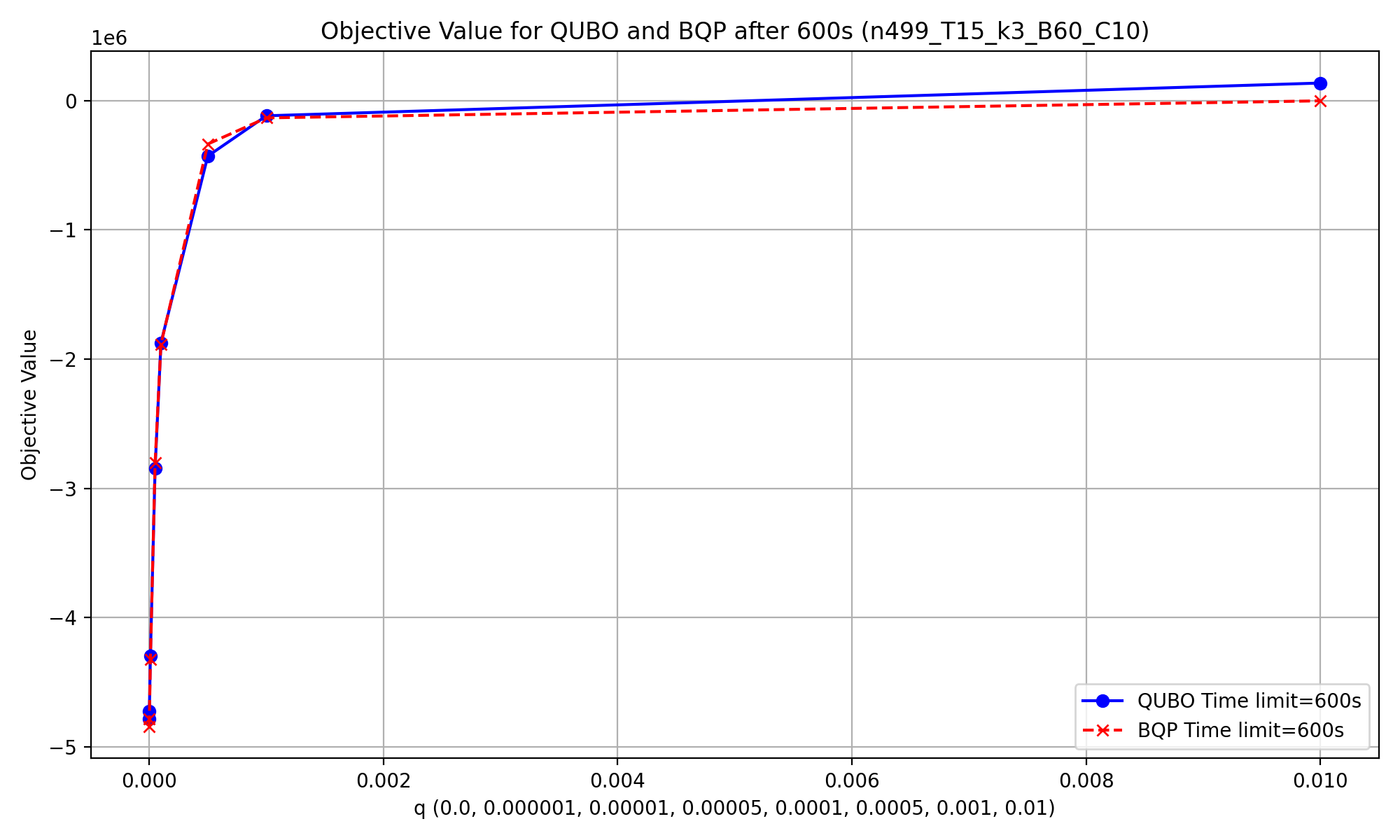}
        \caption{}
        \label{fig:comparison_a499_t15_600s}
    \end{subfigure}
    \vfill
    \begin{subfigure}[b]{0.47\textwidth}
        \centering
        \includegraphics[width=\textwidth]{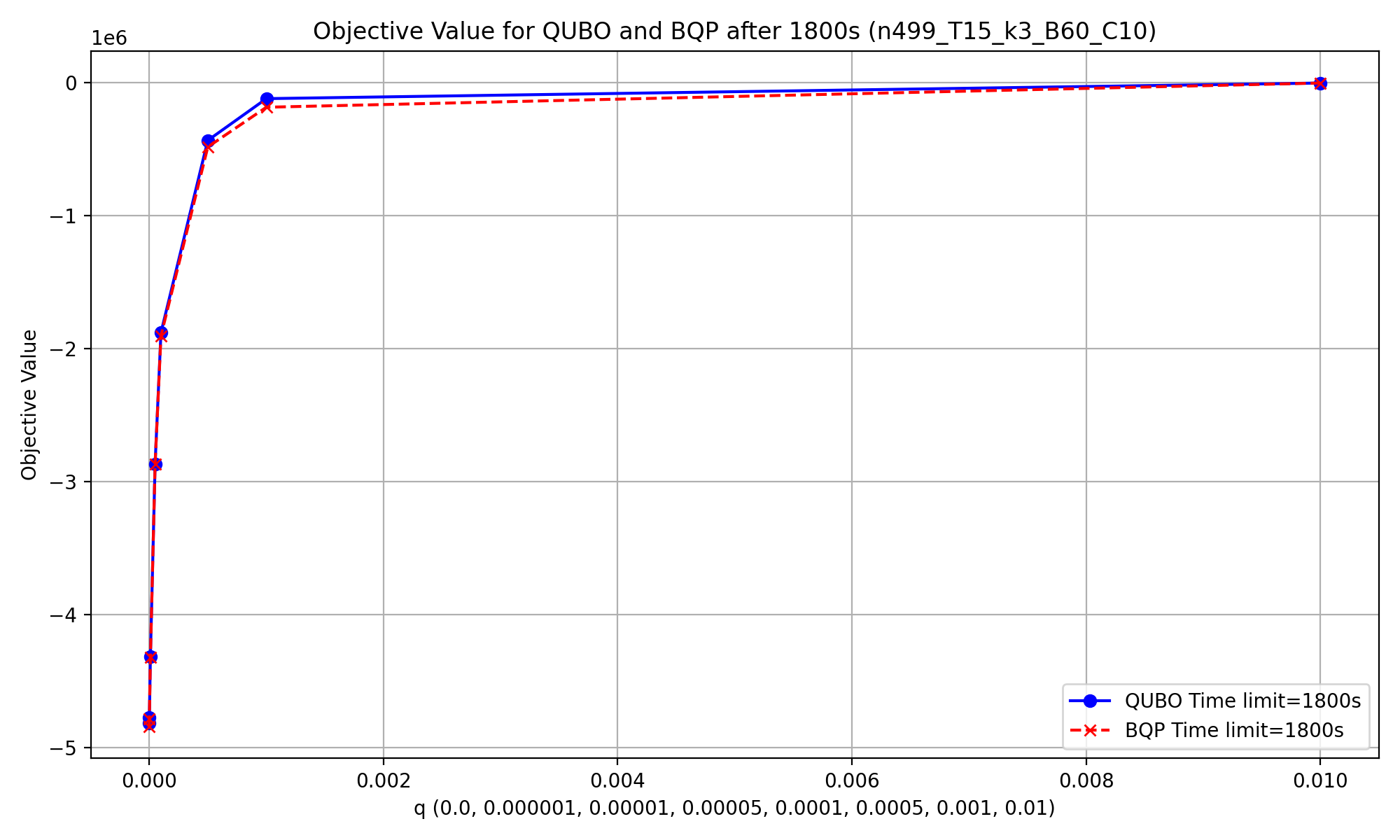}
        \caption{}
        \label{fig:comparison_a499_t15_1800s}
    \end{subfigure}
    \hfill
    \begin{subfigure}[b]{0.47\textwidth}
        \centering
        \includegraphics[width=\textwidth]{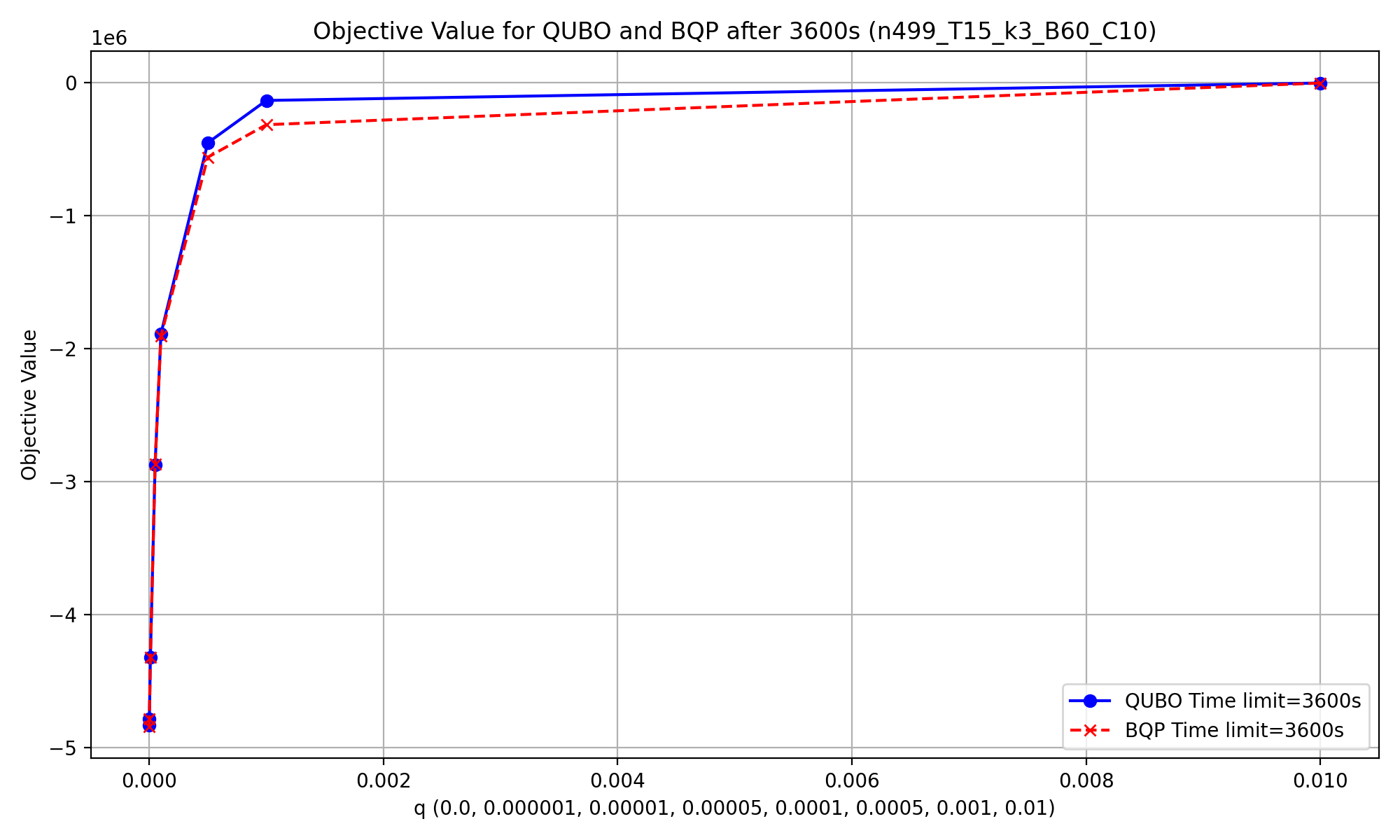}
        \caption{}
        \label{fig:comparison_a499_t15_3600s}
    \end{subfigure}
    \caption{Objective Values for n499\_T15\_k3\_B60\_C10 at different times.}
\end{figure}

\begin{figure}[h!]
    \centering
    \begin{subfigure}[b]{0.47\textwidth}
        \includegraphics[width=\textwidth]{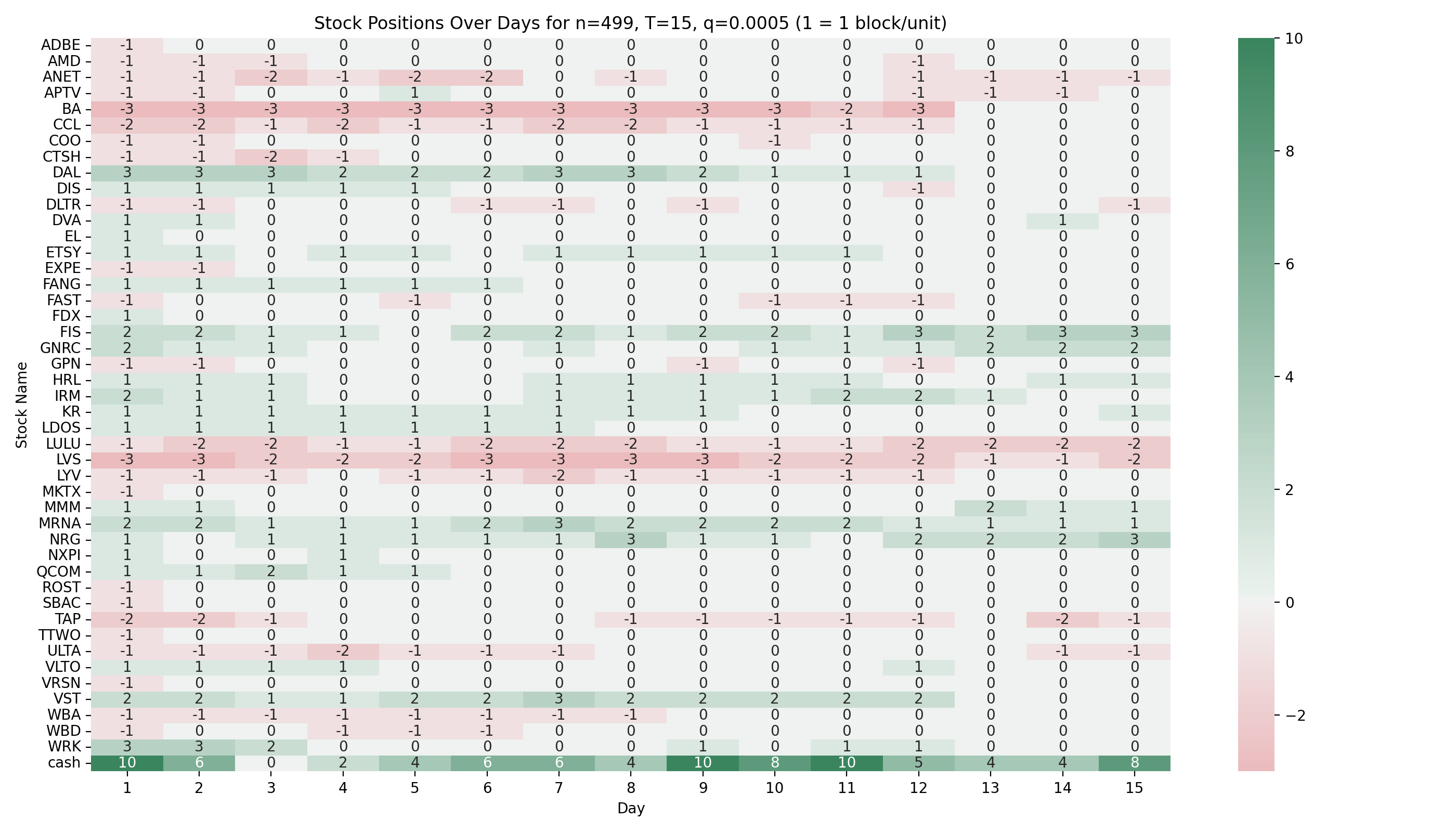}
        \caption{BQP model strategy}
        \label{fig:bqp_499_q0005}
    \end{subfigure}
    \hfill
    \begin{subfigure}[b]{0.47\textwidth}
        \includegraphics[width=\textwidth]{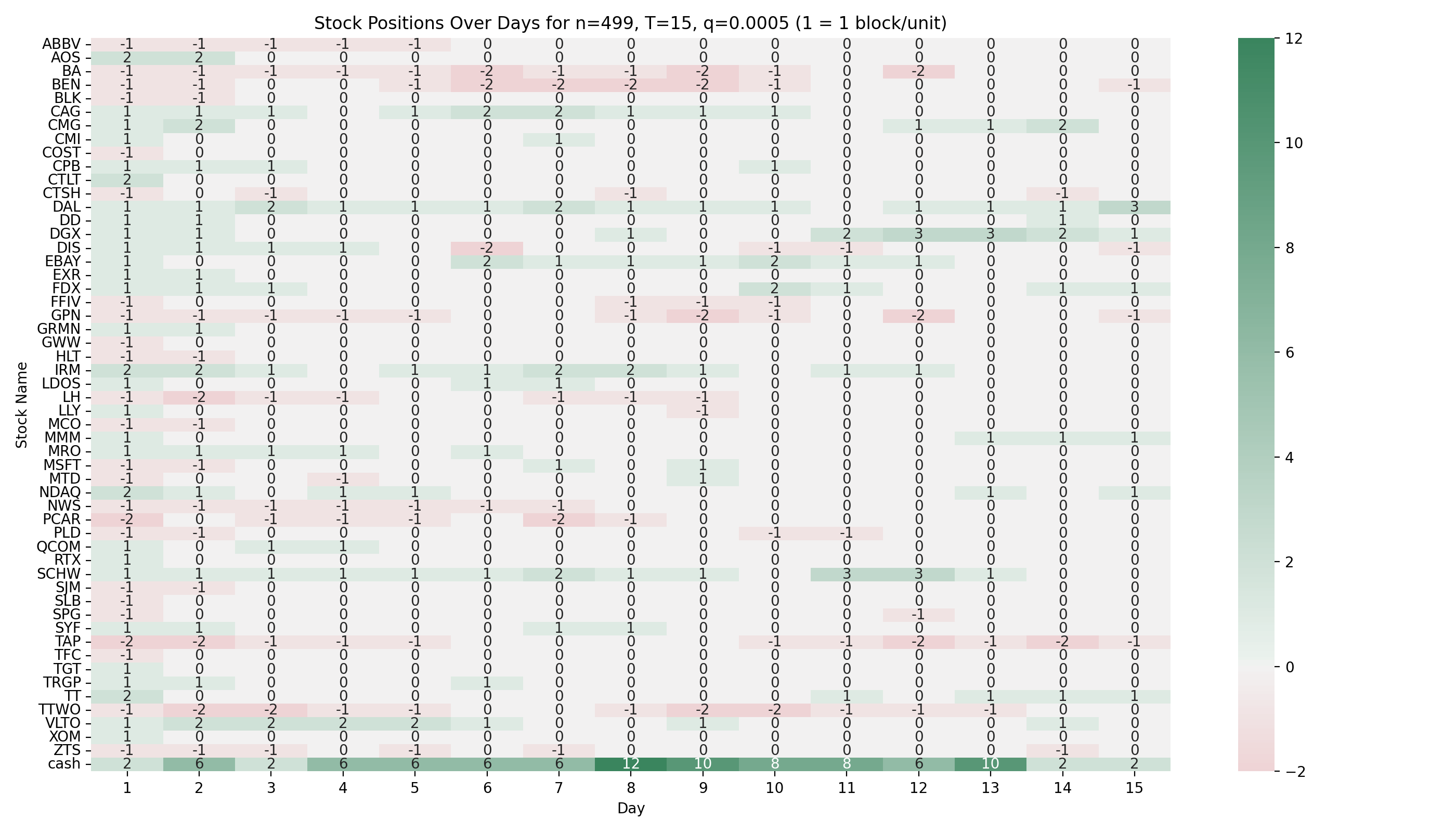}
        \caption{QUBO model strategy}
        \label{fig:qubo_499_q0005}
    \end{subfigure}
    \caption{Stock positions for n499\_T15\_k3\_B60\_C10 with $q = 0.0005$}
\end{figure}

An in-depth analysis of the objective value changes over time for \(q = 0.0005\) and \(q = 0.001\) is depicted in Figures~\ref{fig:a499_t15_q0.0005} and \ref{fig:a499_t15_q0.001}. Although BQP achieves a better final objective value with longer computational time, QUBO reaches good solutions in substantially less time. After 200 seconds for both \(q = 0.0005\) and \(q = 0.001\), there is only marginal improvement, whereas Gurobi requires more than 1400 seconds and 500 seconds, respectively, to stabilize.
We further explored the objective values of QUBO and BQP across different \(q\) values and various computational time limits of 120, 600, 1800, and 3600 seconds, as shown in Figures~\ref{fig:comparison_a499_t15_120s}-\ref{fig:comparison_a499_t15_3600s}. At 120 seconds, except for \(q = 0\), QUBO typically achieves better objective values than BQP, supporting observations that QUBO problems tend to reach better (lower) objective values in the early stages of optimization compared to BQP. At 600 seconds, except for \(q = 0.01\) where the QUBO objective value is noticeably higher, the results for QUBO and BQP are quite similar. By 1800 seconds, the optimization results for QUBO and BQP converge across all \(q\) values.
Figure~\ref{fig:bqp_499_q0005} and Figure~\ref{fig:qubo_499_q0005} illustrate the trading strategy trajectories for the two models. Due to space limitations, we only display the stocks selected on the first day as well as cash positions. The two models suggest differing strategies; the BQP model chose 45 stocks on the first day, while the QUBO model opted for 53 stocks. Nevertheless, their positions on the same stocks throughout the period showed notable similarities. Specifically, both models maintained long positions in DAL from days 1 to 10 and on day 12 and held short positions in BA over the same days, with no positions in BA from days 13 to 15.

The annualized Sharpe ratios and total transaction costs for different \( q \) values of n200\_T10\_k3\_B60\_C10 and n499\_T15\_k3\_B60\_C10 are summarized in Table~\ref{tab:sharpe_ratios}. For \( q \) values less than or equal to 0.0001, the Sharpe ratios for both BQP and QUBO models solved by digital algorithms are very close. In both datasets, Gurobi achieves the highest Sharpe ratio at \( q = 0.001 \), while the ABS2 reaches its peak Sharpe ratio at \( q = 0.0005 \). At \( q = 0.01 \), the high-risk coefficient leads the models to avoid investing in any assets, resulting in a cash-only portfolio with zero variance.

Transaction costs are notably high and similar across models for \( q \) values less than or equal to $0.00001$. As \( q \) increases to $0.0005$ and continues to increase, except for the case of n499\_T15\_k3\_B60\_C10 at \( q = 0.0005 \), the transaction costs for the models decrease significantly with increasing \( q \). This is because as \( q \) increases, the models place greater emphasis on risk, leading to more conservative strategies with lower transaction volumes, and consequently, reduced transaction costs.
\begin{table}[ht]
    \centering
\begin{tabular}{l|c|c|c|c|c|c|c|c}
    \hline\hline
    \hspace{1em} \multirow{3}{*}{\( q \)} & \multicolumn{4}{c|}{n200\_T10\_k3\_B60\_C10} & \multicolumn{4}{c}{n499\_T15\_k3\_B60\_C10}\\
    \cline{2-9}
    & \multicolumn{2}{c|}{Sharpe Ratio} & \multicolumn{2}{c|}{Transaction Cost} & \multicolumn{2}{c|}{Sharpe Ratio} & \multicolumn{2}{c}{Transaction Cost}\\
    \cline{2-9}
    & Gurobi & ABS2 & Gurobi & ABS2 & Gurobi & ABS2 & Gurobi & ABS2 \\
    \hline
    0 & 69.238 & 69.238 & 106,687 & 106,687 & 79.017 & 79.568 & 168,293 & 168,762 \\
    0.000001 & 69.267 & 69.330 & 106,688 & 106,668 & 80.681 & 80.881 & 168,213 & 168,651 \\
    0.00001 & 77.339 & 77.052 & 106,097 & 106,089 & 90.062 & 89.385 & 168,503 & 168,817 \\
    0.00005 & 94.660 & 94.490 & 103,149 & 103,549 & 105.927  &106.228 & 164,742 & 162,685 \\
    0.0001 & 104.813 & 104.637 & 98,608 & 99,242 & 120.484 & 120.256 & 143,231 & 143,464 \\
    0.0005 & 127.370 & 120.346 & 53,628 & 60,486 & 153.603  & 131.520 & 74,240 & 100,364 \\
    0.001 & 135.935 & 103.925 & 42,195 & 30,816 & 165.825 & 100.847 & 63,543 & 125,700 \\
    \hline\hline
\end{tabular}
\caption{Sharpe Ratios and Transaction Costs of n200\_T10\_k3\_B60\_C10 and n499\_T15\_k3\_B60\_C10. When $q=0.01$, the strategy is to hold cash over the whole periods. Thus there is neither variance nor transaction cost.}
\label{tab:sharpe_ratios}
\end{table}
As previously mentioned, models with \( q = 0 \) and \( q = 0.0005 \) were optimized on quantum computers for the dataset n200\_T10\_k3\_B60\_C10. Table~\ref{tab:sharpe_ratios} also compares the economic performance of strategies derived from digital and quantum solvers. The Sharpe ratios obtained from quantum algorithms (D-Wave and InfinityQ) are significantly lower than those achieved by digital algorithms (Gurobi and ABS2), which corresponds to the much higher objective values yielded by D-Wave and InfinityQ compared to Gurobi and ABS2.

In conclusion, while both quantum and digital solvers are capable of effectively handling the portfolio optimization problem, their performance significantly depends on the complexity induced by the risk aversion parameter \(q\). This experiment provides a clear insight into the current computational capabilities and limitations, guiding future developments in optimization algorithms.

\section{Conclusion}\label{sec:conc}
In this study, we have established a robust benchmark for comparing the capabilities of digital and quantum computing in handling complex financial optimization problems. By employing a Quadratic Unconstrained Binary Optimization (QUBO) framework, we addressed dynamic trading strategies for large-scale, multi-period portfolio optimization under realistic market conditions. Our model utilized data from up to 499 stocks in the S\&P 500 over three weeks, factoring in variables such as short selling, transaction costs, cash interest, and integer constraints. This approach enabled us to solve the problem with a precisely determined optimality gap, presenting a standard of solution for the highly intricate models with contemporary computational technologies.

We tested both quantum and digital computing solutions to tackle this optimization challenge. The comparative analysis indicated a negligible performance gap with the recognised digital solver Gurobi, which implies that achieving significantly better results may be unlikely with the current state of technology. Throughout various test scenarios, the qoqo (ABS2) algorithm consistently matched or outperformed Gurobi, highlighting its efficacy. However, improvements in solution quality plateaued after 600 seconds for both methods, with marginal gains.

This benchmarking exercise has demonstrated the QUBO model's effectiveness in portfolio optimization and has highlighted the potential benefits of advancements in quantum computing. Techniques like the Quantum Approximate Optimization Algorithm (QAOA) show potential for leveraging quantum mechanics to solve large-scale multiperiod dynamic optimization. However, current limitations in quantum computing technology, particularly in hardware robustness and noise reduction, pose challenges. As the field advances, these issues are expected to diminish, paving the way for more practical applications of quantum-ready algorithms.

As we release this benchmark, our goal is to enable a clearer understanding of the current state of the art in both quantum and digital computing for financial applications. By doing so, we pave the way for identifying genuine advancements and setting realistic expectations for the performance of emerging computational technologies. Looking forward, as quantum computing matures, it holds the potential to significantly accelerate and enhance the resolution of complex optimization problems, potentially revolutionizing financial modeling and strategy optimization. This advancement will not only improve the efficiency of financial markets but also broaden the scope and applicability of quantum-enhanced solutions in the industry.


\newpage
\clearpage
\printbibliography

@article{lin2008genetic,
  title={Genetic algorithms for portfolio selection problems with minimum transaction lots},
  author={Lin, Chang-Chun and Liu, Yi-Ting},
  journal={European Journal of Operational Research},
  volume={185},
  number={1},
  pages={393--404},
  year={2008},
  publisher={Elsevier}
}

@article{liu2015multi,
  title={A multi-period fuzzy portfolio optimization model with minimum transaction lots},
  author={Liu, Yong-Jun and Zhang, Wei-Guo},
  journal={European Journal of Operational Research},
  volume={242},
  number={3},
  pages={933--941},
  year={2015},
  publisher={Elsevier}
}

@article{gatzert2021portfolio,
  title={Portfolio optimization with irreversible long-term investments in renewable energy under policy risk: A mixed-integer multistage stochastic model and a moving-horizon approach},
  author={Gatzert, Nadine and Martin, Alexander and Schmidt, Martin and Seith, Benjamin and Vogl, Nikolai},
  journal={European Journal of Operational Research},
  volume={290},
  number={2},
  pages={734--748},
  year={2021},
  publisher={Elsevier}
}

@book{mansini2015linear,
  title={Linear and mixed integer programming for portfolio optimization},
  author={Mansini, Renata and WL‚ odzimierz Ogryczak and Speranza, M Grazia and EURO: The Association of European Operational Research Societies},
  volume={21},
  year={2015},
  publisher={Springer}
}

@article{peruzzo2014variational,
  title={A variational eigenvalue solver on a photonic quantum processor},
  author={Peruzzo, Alberto and McClean, Jarrod and Shadbolt, Peter and Yung, Man-Hong and Zhou, Xiao-Qi and Love, Peter J and Aspuru-Guzik, Al{\'a}n and O’brien, Jeremy L},
  journal={Nature communications},
  volume={5},
  number={1},
  pages={4213},
  year={2014},
  publisher={Nature Publishing Group UK London}
}

@article{braine2021quantum,
  title={Quantum algorithms for mixed binary optimization applied to transaction settlement},
  author={Braine, Lee and Egger, Daniel J and Glick, Jennifer and Woerner, Stefan},
  journal={IEEE Transactions on Quantum Engineering},
  volume={2},
  pages={1--8},
  year={2021},
  publisher={IEEE}
}

@article{gilliam2021grover,
  title={Grover adaptive search for constrained polynomial binary optimization},
  author={Gilliam, Austin and Woerner, Stefan and Gonciulea, Constantin},
  journal={Quantum},
  volume={5},
  pages={428},
  year={2021},
  publisher={Verein zur F{\"o}rderung des Open Access Publizierens in den Quantenwissenschaften}
}

@article{herman2022survey,
  title={A survey of quantum computing for finance},
  author={Herman, Dylan and Googin, Cody and Liu, Xiaoyuan and Galda, Alexey and Safro, Ilya and Sun, Yue and Pistoia, Marco and Alexeev, Yuri},
  journal={arXiv preprint arXiv:2201.02773},
  year={2022}
}

@article{mugel2022dynamic,
  title={Dynamic portfolio optimization with real datasets using quantum processors and quantum-inspired tensor networks},
  author={Mugel, Samuel and Kuchkovsky, Carlos and Sanchez, Escolastico and Fernandez-Lorenzo, Samuel and Luis-Hita, Jorge and Lizaso, Enrique and Orus, Roman},
  journal={Physical Review Research},
  volume={4},
  number={1},
  pages={013006},
  year={2022},
  publisher={APS}
}

@techreport{elsokkary2017financial,
  title={Financial portfolio management using d-wave quantum optimizer: The case of abu dhabi securities exchange},
  author={Elsokkary, Nada and Khan, Faisal Shah and La Torre, Davide and Humble, Travis S and Gottlieb, Joel},
  year={2017},
  institution={Oak Ridge National Lab.(ORNL), Oak Ridge, TN (United States)}
}

@article{lucas2014ising,
  title={Ising formulations of many NP problems},
  author={Lucas, Andrew},
  journal={Frontiers in physics},
  volume={2},
  pages={5},
  year={2014},
  publisher={Frontiers}
}

@article{bertsimas1998optimal,
  title={Optimal control of execution costs},
  author={Bertsimas, Dimitris and Lo, Andrew W},
  journal={Journal of financial markets},
  volume={1},
  number={1},
  pages={1--50},
  year={1998},
  publisher={Elsevier}
}

@book{bertsimas1997introduction,
  title={Introduction to linear optimization},
  author={Bertsimas, Dimitris and Tsitsiklis, John N},
  volume={6},
  year={1997},
  publisher={Athena scientific Belmont, MA}
}

@article{brandt2006dynamic,
  title={Dynamic portfolio selection by augmenting the asset space},
  author={Brandt, Michael W and Santa-Clara, Pedro},
  journal={The journal of Finance},
  volume={61},
  number={5},
  pages={2187--2217},
  year={2006},
  publisher={Wiley Online Library}
}

@article{broadie1998optimal,
  title={Optimal replication of contingent claims under portfolio constraints},
  author={Broadie, Mark and Cvitani{\'c}, Jak{\v{s}}a and Soner, H Mete},
  journal={The Review of Financial Studies},
  volume={11},
  number={1},
  pages={59--79},
  year={1998},
  publisher={Oxford University Press}
}

@article{campbell1998econometrics,
  title={The econometrics of financial markets},
  author={Campbell, John Y and Lo, Andrew W and MacKinlay, A Craig and Whitelaw, Robert F},
  journal={Macroeconomic Dynamics},
  volume={2},
  number={4},
  pages={559--562},
  year={1998},
  publisher={Cambridge University Press}
}

@article{dybvig1985differential,
  title={Differential information and performance measurement using a security market line},
  author={Dybvig, Philip H and Ross, Stephen A},
  journal={The Journal of finance},
  volume={40},
  number={2},
  pages={383--399},
  year={1985},
  publisher={Wiley Online Library}
}

@article{egger2020credit,
  title={Credit risk analysis using quantum computers},
  author={Egger, Daniel J and Guti{\'e}rrez, Ricardo Garc{\'\i}a and Mestre, Jordi Cahu{\'e} and Woerner, Stefan},
  journal={IEEE Transactions on Computers},
  volume={70},
  number={12},
  pages={2136--2145},
  year={2020},
  publisher={IEEE}
}

@article{fama1970efficient,
  title={Efficient capital markets: A review of theory and empirical work},
  author={Fama, Eugene F},
  journal={The journal of Finance},
  volume={25},
  number={2},
  pages={383--417},
  year={1970},
  publisher={JSTOR}
}

@article{harrow2017quantum,
  title={Quantum computational supremacy},
  author={Harrow, Aram W and Montanaro, Ashley},
  journal={Nature},
  volume={549},
  number={7671},
  pages={203--209},
  year={2017},
  publisher={Nature Publishing Group UK London}
}

@article{li2000optimal,
  title={Optimal dynamic portfolio selection: Multiperiod mean-variance formulation},
  author={Li, Duan and Ng, Wan-Lung},
  journal={Mathematical finance},
  volume={10},
  number={3},
  pages={387--406},
  year={2000},
  publisher={Wiley Online Library}
}

@incollection{merton1975optimum,
  title={Optimum consumption and portfolio rules in a continuous-time model},
  author={Merton, Robert C},
  booktitle={Stochastic optimization models in finance},
  pages={621--661},
  year={1975},
  publisher={Elsevier}
}

@article{orus2019quantum,
  title={Quantum computing for finance: Overview and prospects},
  author={Or{\'u}s, Rom{\'a}n and Mugel, Samuel and Lizaso, Enrique},
  journal={Reviews in Physics},
  volume={4},
  pages={100028},
  year={2019},
  publisher={Elsevier}
}

@book{pliska1997introduction,
  title={Introduction to mathematical finance: discrete time models},
  author={Pliska, Stanley R},
  year={1997},
  publisher={Wiley}
}

@article{rebentrost2018quantum,
  title={Quantum computational finance: Monte Carlo pricing of financial derivatives},
  author={Rebentrost, Patrick and Gupt, Brajesh and Bromley, Thomas R},
  journal={Physical Review A},
  volume={98},
  number={2},
  pages={022321},
  year={2018},
  publisher={APS}
}

@article{samuelson1975lifetime,
  title={Lifetime portfolio selection by dynamic stochastic programming},
  author={Samuelson, Paul A},
  journal={Stochastic optimization models in finance},
  pages={517--524},
  year={1975},
  publisher={Elsevier}
}

@article{sharpe1998sharpe,
  title={The sharpe ratio},
  author={Sharpe, William F},
  journal={Streetwise--the Best of the Journal of Portfolio Management},
  volume={3},
  pages={169--85},
  year={1998},
  publisher={Princeton University Press NJ}
}

@article{stamatopoulos2020option,
  title={Option pricing using quantum computers},
  author={Stamatopoulos, Nikitas and Egger, Daniel J and Sun, Yue and Zoufal, Christa and Iten, Raban and Shen, Ning and Woerner, Stefan},
  journal={Quantum},
  volume={4},
  pages={291},
  year={2020},
  publisher={Verein zur F{\"o}rderung des Open Access Publizierens in den Quantenwissenschaften}
}

@article{woerner2019quantum,
  title={Quantum risk analysis},
  author={Wörner, Stefan and Egger, Daniel J},
  journal={npj Quantum Information},
  volume={5},
  number={1},
  pages={15},
  year={2019},
  publisher={Nature Publishing Group UK London}
}

@article{zabarankin2014capital,
  title={Capital asset pricing model (CAPM) with drawdown measure},
  author={Zabarankin, Michael and Pavlikov, Konstantin and Uryasev, Stan},
  journal={European Journal of Operational Research},
  volume={234},
  number={2},
  pages={508--517},
  year={2014},
  publisher={Elsevier}
}

@article{dai2021learning,
  title={Learning equilibrium mean-variance strategy},
  author={Dai, Min and Dong, Yuchao and Jia, Yanwei},
  journal={Mathematical Finance},
  year={2021},
  publisher={Wiley Online Library}
}

@article{brown2022information,
  title={Information relaxations and duality in stochastic dynamic programs: A review and tutorial},
  author={Brown, David B and Smith, James E and others},
  journal={Foundations and Trends{\textregistered} in Optimization},
  volume={5},
  number={3},
  pages={246--339},
  year={2022},
  publisher={Now Publishers, Inc.}
}

@article{brown2011dynamic,
  title={Dynamic portfolio optimization with transaction costs: Heuristics and dual bounds},
  author={Brown, David B and Smith, James E},
  journal={Management Science},
  volume={57},
  number={10},
  pages={1752--1770},
  year={2011},
  publisher={INFORMS}
}

@book{nielsen2010quantum,
  title={Quantum computation and quantum information},
  author={Nielsen, Michael A and Chuang, Isaac L},
  year={2010},
  publisher={Cambridge university press}
}

@book{mcgeoch2022adiabatic,
  title={Adiabatic quantum computation and quantum annealing: Theory and practice},
  author={McGeoch, Catherine C},
  year={2022},
  publisher={Springer Nature}
}

@inproceedings{kerenidis2019quantum,
  title={Quantum algorithms for portfolio optimization},
  author={Kerenidis, Iordanis and Prakash, Anupam and Szil{\'a}gyi, D{\'a}niel},
  booktitle={Proceedings of the 1st ACM Conference on Advances in Financial Technologies},
  pages={147--155},
  year={2019}
}

@ARTICLE{rosenberg2016solving,
  author={Rosenberg, Gili and Haghnegahdar, Poya and Goddard, Phil and Carr, Peter and Wu, Kesheng and de Prado, Marcos López},
  journal={IEEE Journal of Selected Topics in Signal Processing}, 
  title={Solving the Optimal Trading Trajectory Problem Using a Quantum Annealer}, 
  year={2016},
  volume={10},
  number={6},
  pages={1053-1060},
  doi={10.1109/JSTSP.2016.2574703}}

@article{gaivoronski2005optimal,
  title={Optimal portfolio selection and dynamic benchmark tracking},
  author={Gaivoronski, Alexei A and Krylov, Sergiy and Van der Wijst, Nico},
  journal={European Journal of operational research},
  volume={163},
  number={1},
  pages={115--131},
  year={2005},
  publisher={Elsevier}
}

@article{rietz1988equity,
  title={The equity risk premium a solution},
  author={Rietz, Thomas A},
  journal={Journal of monetary Economics},
  volume={22},
  number={1},
  pages={117--131},
  year={1988},
  publisher={Elsevier}
}

@article{d2002market,
  title={The market for borrowing stock},
  author={D’avolio, Gene},
  journal={Journal of financial economics},
  volume={66},
  number={2-3},
  pages={271--306},
  year={2002},
  publisher={Elsevier}
}

@article{tepla2000optimal,
  title={Optimal portfolio policies with borrowing and shortsale constraints},
  author={Tepla, Lucie},
  journal={Journal of Economic Dynamics and Control},
  volume={24},
  number={11-12},
  pages={1623--1639},
  year={2000},
  publisher={Elsevier}
}

@article{rehfeldt2023faster,
  title={Faster exact solution of sparse MaxCut and QUBO problems},
  author={Rehfeldt, Daniel and Koch, Thorsten and Shinano, Yuji},
  journal={Mathematical Programming Computation},
  volume={15},
  number={3},
  pages={445--470},
  year={2023},
  publisher={Springer},
  doi={10.1007/s12532-023-00236-6}
}

@inproceedings{nakano2023diverse,
  title={Diverse Adaptive Bulk Search: a Framework for Solving QUBO Problems on Multiple GPUs},
  author={Nakano, Katsuki and others},
  booktitle={2023 IEEE International Parallel and Distributed Processing Symposium Workshops (IPDPSW)},
  pages={393--402},
  year={2023},
  organization={IEEE},
  doi={10.1109/IPDPSW59300.2023.00060}
}

@article{montanez2024towards,
  title={Towards a universal QAOA protocol: Evidence of quantum advantage in solving combinatorial optimization problems},
  author={Montanez-Barrera, JA and Michielsen, Kristel},
  journal={arXiv preprint arXiv:2405.09169},
  year={2024}
}

\end{document}